\documentclass[final,10pt]{elsarticle}
\usepackage{algorithm,algorithmic,appendix}
\usepackage{cite}
\usepackage{fullpage}
\usepackage{multirow}
\usepackage[latin1]{inputenc}
\usepackage[crop=pdfcrop]{pstool}

\usepackage{color}
\newcommand{\state}{\ensuremath{{\boldsymbol x}}}	% 
\newcommand{\timestep}{\ensuremath{{h}}}	% 
\newcommand{\snapsNo}{\ensuremath{\mathcal X}}	% 
\newcommand{\snapmat} {\ensuremath{\boldsymbol W}}
\newcommand{\energyCrit} {\ensuremath{\nu}}
\newcommand{\snaps}[1]{\ensuremath{\snapsNo_{#1}}}	% 
\newcommand{\podArgsNo}{\ensuremath{\boldsymbol \Phi^e}}	% 
\newcommand{\podArgs}[2]{\ensuremath{\podArgsNo\left(#1,#2\right)}}	% 
\newcommand{\basisIndex}{\ensuremath{j}}	% 
\newcommand{\nTrain}{\ensuremath{n_{\text{train}}}}	% 
\newcommand{\efomone}{\ensuremath{\varepsilon}}	% 
\newcommand{\FOM}{\ensuremath{{\text{FOM}}}}	% 
\newcommand{\walltime}{\ensuremath{\text{T}}}	% 
\newcommand{\temporalRedFactor}{\ensuremath{\kappa}}	% 
\newcommand{\temporalRedFactorImprove}{\ensuremath{\mathfrak k}}	% 
\newcommand{\speedup}{\ensuremath{S}}	% 
\newcommand{\speedupImprove}{\ensuremath{\mathfrak s}}	% 
\newcommand{\Dtrain}{\ensuremath{\paramDomain_{\text{train}}}}	% 
\newcommand{\paramTrain}{\Dtrain}	% 
\newcommand{\snapsMat}[1]{\ensuremath{ \boldsymbol Y_{#1}}}	% 
\newcommand{\redSnapsMat}[1]{\ensuremath{ \hat {\boldsymbol Y}_{#1}}}	% 
\newcommand{\redSnapsVec}[2]{\ensuremath{ \hat {\boldsymbol y}_{#2,#1}}}	% 

\newcommand{\ones}{\ensuremath{\mathbf 1}}
\newcommand{\avgNew}{\ensuremath{ K}}	% 
\newcommand{\restrict}{\ensuremath{\boldsymbol Z}}	% 
\newcommand{\restrictn}{\ensuremath{\boldsymbol Z(n,\memory)}}	% 
\newcommand{\Rapprox}{\ensuremath{\tilde {\boldsymbol r}}}	% 
\newcommand{\totalTS}{\ensuremath{M}}	% 
\newcommand{\totalTSun}{\ensuremath{M}}	% 
\newcommand{\totalLinSys}{\ensuremath{\avgNew\totalTS}}	% 
\newcommand{\sparsity}{\ensuremath{\omega}}	% 
\newcommand{\linearIt}{\ensuremath{ L}}	% 
\newcommand{\timebasis}{\ensuremath{\boldsymbol \Xi}}	% 
\newcommand{\timebasisj}{\ensuremath{\timebasis_\basisIndex}}	% 
\newcommand{\unrollNone}{\ensuremath{\underline h}}	% 
\newcommand{\unroll}[1]{\ensuremath{\unrollNone\left(#1\right)}}	% 
\newcommand{\unrollinv}[1]{\ensuremath{\unrollNone^{-1}\left(#1\right)}}	% 
\newcommand{\approxstate}{\ensuremath{\tilde{\state}}}	% 
\newcommand{\forcing}{\ensuremath{\boldsymbol p}}	% 
\newcommand{\nforcing}{\ensuremath{\mathrm p}}	% 
\newcommand{\forceDist}{\ensuremath{{\mathbf r}}}	% 
\newcommand{\forceTime}{\ensuremath{{p}}}	% 
\newcommand{\unknown}{\ensuremath{\boldsymbol w}}	% 
\newcommand{\unknownTime}{\ensuremath{\unknown}}	% 
\newcommand{\unknownTimeInitk}{\ensuremath{\unknownRefArg{\param_k}}}	% 
\newcommand{\unknownApprox}{\ensuremath{\tilde\unknown}}	% 
\newcommand{\dummy}{\ensuremath{\boldsymbol y}}	% 
\newcommand{\param}{\ensuremath{\boldsymbol q}}	% 
\newcommand{\parami}[1]{\ensuremath{ q_{#1}}}	% 
\newcommand{\nstate}{\ensuremath{{\hat N}}}	% 
\newcommand{\ntime}{\ensuremath{a}}	% 
\newcommand{\ntimej}{\ensuremath{{\ntime_j}}}	% 
\newcommand{\podstate}{\ensuremath{{\boldsymbol \Phi}}}	% 
\newcommand{\testBasis}{\ensuremath{{\boldsymbol \Psi}}}	% 
\newcommand{\podstatevec}{\ensuremath{\boldsymbol \phi}}	% 
\newcommand{\redstate}{\ensuremath{{\hat {\state}}}}	% 
\newcommand{\redstateEntry}[1]{\ensuremath{{\hat {x}_{#1}}}}	% 
\newcommand{\redtime}{\ensuremath{\boldsymbol z}}	% 
\newcommand{\redtimej}{\ensuremath{\redtime_\basisIndex}}	% 
\newcommand{\podres}{\ensuremath{\boldsymbol \Phi_R}}	% 
\newcommand{\podjac}{\ensuremath{\boldsymbol \Phi_J}}	% 
\newcommand{\podf}{\ensuremath{\boldsymbol \Phi_f}}	% 
\newcommand{\podg}{\ensuremath{\boldsymbol \Phi_g}}	% 
\newcommand{\nt}{\ensuremath{\totalTS}}	% 
\newcommand{\memory}{\ensuremath{\alpha}}	% 
\newcommand{\memoryMax}{\ensuremath{\alpha_{\max}}}	% 
\newcommand{\newtonThreshold}{\ensuremath{\tau}}	% 
\newcommand{\newtonIts}{\ensuremath{K}}	% 
\newcommand{\newtonItsAvg}{\ensuremath{\bar \newtonIts}}	% 
\newcommand{\forecastVar}{\ensuremath{\underline \redUnknown}}	% 
\newcommand{\forecastTime}{\ensuremath{ \forecastVar}}	% 
\newcommand{\forecastTimej}{\ensuremath{\underline {\hat w}_j}}	% 
\newcommand{\nparam}{\ensuremath{\mathrm q}}

\newcommand{\lastT}{\ensuremath{ T}}	% 

\newcommand{\sampleMat} {\ensuremath{\boldsymbol Z}}
\newcommand{\nsample} {\ensuremath{n_Z}}

\newcommand{\paramOnline}{\ensuremath{\param^\star}}	% 
\newcommand{\paramOnlineNum}[1]{\ensuremath{\param^{\star,#1}}}	% 
\newcommand{\nenergy}{\ensuremath{n_e}}
\newcommand{\nsnap} {\ensuremath{{n_w}}}

\newcommand{\stateInitialNo} {\ensuremath{\state^0}}
\newcommand{\stateInitial} {\ensuremath{\state^0\left(\param\right)}}
\newcommand{\stateRefNo} {\ensuremath{\bar\state}}
\newcommand{\stateRefArg}[1]{\ensuremath{\stateRefNo\left(#1\right)}}
\newcommand{\stateRef} {\ensuremath{\stateRefArg{\param}}}
\newcommand{\unknownRefNo} {\ensuremath{\bar\unknown}}
\newcommand{\unknownRefArg}[1]{\ensuremath{\unknownRefNo\left(#1\right)}}
\newcommand{\unknownRef} {\ensuremath{\unknownRefArg{\param}}}
\newcommand{\residual} {\ensuremath{\boldsymbol r}}
\newcommand{\Rn} {\residual^n}
\newcommand{\tildeRn} {\tilde\residual^n}
\newcommand{\Rni} {\ensuremath{\residual^n_i}}

\newcommand{\unknownnarg}[1] {\ensuremath{\unknown^{n,#1}}}
\newcommand{\unknownnone} {\ensuremath{\unknownnarg{1}}}
\newcommand{\unknownns} {\ensuremath{\unknownnarg{s}}}
\newcommand{\unknownni} {\ensuremath{\unknownnarg{i}}}
\newcommand{\redUnknown}{\ensuremath{ \hat \unknown}}	% 

\newcommand{\redUnknownn} {\ensuremath{\redUnknown^n}}

\newcommand{\redUnknownTime}{\ensuremath{ \hat \unknownTime}}	% 
\newcommand{\redUnknownTimeEntry}[1]{\ensuremath{ \hat w_{#1}}}	% 
\newcommand{\redUnknownTimeZeroEntry}[1]{\ensuremath{ \hat w_{#1}^{n(0)}}}	% 
\newcommand{\rayleighMass}{\ensuremath{\alpha}}	% 
\newcommand{\rayleighStiff}{\ensuremath{\beta}}	% 

\newcommand{\Mass} {\ensuremath{{\boldsymbol M}}}
\newcommand{\Mparam} {\ensuremath{{\Mass\left(\param\right)}}}

\newcommand{\Conly} {\ensuremath{\boldsymbol C}}
\newcommand{\Cparam} {\ensuremath{{\Conly\left(\param\right)}}}

\newcommand{\fext} {\ensuremath{\boldsymbol f_\text{ext}}}

\newcommand{\fint} {\ensuremath{\boldsymbol f_\text{int}}}
\newcommand{\fintparam} {\ensuremath{\fint\left(\state;\param\right)}}
\newcommand{\fintparamArg}[1]{\ensuremath{\fint\left(#1;\param\right)}}
\newcommand{\fintparamArgs}[2]{\ensuremath{\fint\left(#1;#2\right)}}

\newcommand{\paramDomain}{\ensuremath{\mathcal D}}
\newcommand{\forceMag}{\ensuremath{\gamma}}

\newcommand{\f} {\ensuremath{{\boldsymbol f}}}
\newcommand{\g} {\ensuremath{{\boldsymbol g}}}

\newcommand{\initialCond}{\ensuremath{{\boldsymbol s}}}	% 
\newcommand{\initialCondMag}{\ensuremath{{ s}}}	% 
\newcommand{\paramNom}{\ensuremath{\bar\param}}
\newcommand{\fNomi}{\ensuremath{\underline f_i}}	% 
\newcommand{\fNomnum}[1]{\ensuremath{\underline f_{#1}}}	% 

\newcommand{\vVec} {\ensuremath{{\boldsymbol v}}}
\newcommand{\w} {\ensuremath{{\boldsymbol w}}}

\newcommand{\M} {\ensuremath{{\boldsymbol M}}}
\newcommand{\identity} {\ensuremath{{\boldsymbol I}}}
\newcommand{\ndof} {\ensuremath{{N}}}
\newcommand{\vInit} {\ensuremath{{\vVec^0}}}

\newcounter{lemctr}
\newcounter{pfctr}
\newcounter{remctr}
\newenvironment{remark}{\refstepcounter{remctr}\begin{trivlist}
\item {\emph {Remark}.\:}}
{\end{trivlist}}

\newcounter{propctr}
\newcounter{proposctr}
\newcommand{\RR}[1]{\ensuremath{\mathbb{R}^{ #1 }}}

\newcommand{\vecmat}[2]{\left[{{#1}}_1 \ \cdots\ {{#1}}_{#2}\right]}
\newcommand{\vecmatsub}[3]{\left[{{#1}_{#2,1} \ \cdots\ {#1}_{#2,#3}}\right]}

\usepackage{algorithmic, algorithm}% subcaptions for subfigures
\usepackage{amsmath,amssymb}% subcaptions for subfigures
\usepackage{graphics,epsfig}% subcaptions for subfigures
\usepackage{subfigmat}% matrices of similar subfigures, aka small mulitples
\usepackage{amsopn}

\usepackage[hmargin=3cm,vmargin=3.5cm]{geometry}
\usepackage{hyperref}

\begin{document}

\begin{frontmatter}
\title{Decreasing the temporal complexity for\\ nonlinear, implicit
reduced-order models by forecasting }
\author{Kevin Carlberg, Jaideep Ray, and Bart van Bloemen Waanders \\ Sandia
National Laboratories\footnote{Sandia is a multiprogram laboratory operated by
Sandia Corporation, a Lockheed Martin Company, for the United States
Department of Energy under contract DE-AC04-94-AL85000.}\\
7011 East Ave, MS 9159,  Livermore, CA 94550}

%\date{}

\begin{abstract}

\noindent Implicit numerical integration of nonlinear ODEs requires solving a
system of nonlinear algebraic equations at each time step.  Each of these systems is
often solved by a Newton-like method, which incurs a sequence of linear-system
solves. Most model-reduction techniques for nonlinear ODEs exploit knowledge
of system's spatial behavior to reduce the computational complexity of each
linear-system solve. However, the number of linear-system solves for the
reduced-order simulation often remains roughly the same as that for the
full-order simulation.	

We propose exploiting knowledge of the model's temporal behavior to 1)
forecast the unknown variable of the reduced-order system of nonlinear
equations at future time steps, and 2) use this forecast as an initial guess
for the Newton-like solver during the reduced-order-model simulation. To compute
the forecast, we propose using the Gappy POD technique. The goal is to
generate an accurate initial guess so that the Newton solver requires many
fewer iterations to converge, thereby decreasing the number of
linear-system solves in the reduced-order-model simulation.
\end{abstract}


	\end{frontmatter}

\section{Introduction}\label{sec:intro}

%\subsection{Motivation}

High-fidelity physics-based numerical simulation has become an indispensable
engineering tool across a wide range of disciplines. Unfortunately, such
simulations often bear an extremely large computational cost due to the
large-scale, nonlinear nature of many high-fidelity models. When an implicit 
integrator is employed to advance the solution in time (as is often essential,
e.g., for stiff problems) this large cost arises from the need to solve a
sequence of high-dimensional systems of nonlinear algebraic equations---one at
each time step. As a result, individual simulations can take weeks or months
to complete, even when high-performance computing resources are available.
This renders such simulations impractical for time-critical and many-query
applications.  For example, uncertainty-quantification applications (e.g.,
Bayesian inference problems) call for hundreds or thousands of simulations
(i.e., forward solves) to be completed in days or weeks; in-the-field analysis
(e.g., guidance in-field data acquisition) requires near-real-time simulation.

Projection-based nonlinear model-reduction techniques have been successfully
applied to decrease the computational cost of high-fidelity simulation while
retaining high levels of accuracy. To accomplish this, these methods exploit
knowledge of the system's dominant \emph{spatial behavior}---as observed during
`training simulations' conducted \emph{a priori}---to decrease the simulation's
\emph{spatial complexity}, which we define as the computational cost of each
linear-system solve.\footnote{A sequence of linear systems arises at each time
step when a Newton-like method is employed to solve the system of nonlinear
algebraic equations.} To do so, these methods 1) decrease the dimensionality of the
linear systems by projection, and 2) approximate vector-valued nonlinear functions by
sampling methods that compute only a few of the vector's entries (e.g.,
empirical interpolation \citep{barrault2004eim,chaturantabut2010journal},
Gappy POD \citep{sirovichOrigGappy}).  However, these techniques are often
insufficient to adequately reduce the computational cost of the simulation.
For example, Ref. \citep{carlbergJCP} presented results for the GNAT nonlinear
model-reduction technique applied to a large-scale nonlinear turbulent-flow
problem. The reduced-order model generated solutions with sub-1\% errors,
reduced the spatial complexity by a factor of 637, and employed only 4
computing cores---a significant reduction from the 512 cores required for the
high-fidelity simulation. However, the total number of
linear-system solves required for the reduced-order-model simulation, which we
define as the \emph{temporal complexity}, remained large. In fact, the
temporal complexity was decreased by a factor of only 1.5.  As a result, the
total computing resources (computing cores $\times$ wall time) required for
the simulation were decreased by a factor of 438, but the wall time was
reduced by a factor of merely 6.9. While these results are promising (especially
in their ability to reduce spatial complexity), the time integration of
nonlinear dynamics remains problematic and often precludes real-time
performance.

The goal of this work is exploit knowledge of the system's \emph{temporal
behavior} as observed during the training simulations to decrease the temporal
complexity of reduced-order-model simulations. For this
purpose, we first briefly review methods that exploit observed temporal
behavior to improve computational performance.

Temporal forecasting techniques have been investigated for many years with a
specific focus on reducing wall time in a stable manner with maximal accuracy.
The associated body of work is large and a comprehensive review is beyond the
scope of this paper. However, this work focuses on time integration for
reduced-order models of highly nonlinear dynamical systems; several categories
of specialized research efforts provide an appropriate context for this
research.

At the most fundamental level of temporal forecasting, a variety of
statistical time-series-analysis methods exist that exploit 1) knowledge of
the temporal structure, e.g., smoothness, of a model's  variables, and 2)
previous values of these variables for the current time series or trajectory.
The connection between these methods and our work is that such forecasts can
serve as an initial guess for an iterative solver (e.g., Newton's method) at
an advanced point in time. However, the disconnect between such methods and
the present context is that randomness and uncertainty drive time-series
analysis; as such, these forecasting methods are stochastic in nature (see
Refs.\ \citep{90dp1a, 97gc1a,82er1a,61gf1a,73hw2a,arima, Holt20045,60wp1a}).
In addition, the majority of time-series analyses have been applied to
application domains (e.g., economics) with dynamics that are not generally
modeled using partial differential equations. Finally, such forecasting
techniques do not exploit a collection of observed, complete time histories
from training experiments conducted \emph{a priori}.  Because such training
simulations lend important insight into the spatial and temporal behavior of
the model, we are interested in developing a technique that can exploit such
data.

Alternatively, time integrators for ordinary
differential equations (ODEs) employ polynomial extrapolations to
provide reasonably accurate forecasts of the state or the unknown at
each time step. Time integrators employ such a forecast for two
purposes.  First, algorithms with adaptive time steps employ
interpolation to obtain solutions (and their time derivatives) at
arbitrary points in time. Implicit time integrators for nonlinear
ODEs, which require the iterative solution of nonlinear algebraic
systems at each time step, use recent history (of the current
trajectory) to forecast an accurate guess of the unknown in the
algebraic system (see, e.g., Ref.~\citep{doi:10.1137/0910062}).
Again, forecasting by polynomial extrapolation makes no use of the temporal
behavior observed during training simulations.

Closely connected to time integration but specialized to leverage developments
in high-performance computing, time parallel methods can offer computational
speedup when integrating ODEs.  Dating back to before the general availability
of parallel computers, researchers speculated about the benefits of decomposing
the temporal domain across multiple processors \citep{Nievergelt_ACM_64}.
Advancements have been made from parallel multigrid to parareal techniques
\citep{Gander_NLAA_93,Horton_SISC_95,lions2001parareal,cortial2011time}.  Although
time-domain decomposition algorithms have demonstrated speedup, they are
limited in comparison to the spatial domain decomposition methods and they
require a careful balance between stability and computational efficiency
\citep{farhat2003time}.  It is possible that these methods could further
improve performance in a model-reduction setting \citep{harden2008real}
(and could complement the method proposed in this work), but near real-time
performance is likely unachievable through time-parallel methods alone.

To some extent, exploiting temporal behavior has been explored in nonlinear
model reduction.  Bos et al.\ \citep{bos2004als} proposed a reduced-order model
in the context of explicit time integration wherein the generalized coordinates
are computed based on a best-linear-unbiased (BLU) estimate approach. Here, the
reduced state coordinates at time step $n+1$ are computed using empirically
derived correlations between the reduced state coordinates and 1) their value
at the previous time step, 2) the forcing input at the previous time step, and
3) a subset of the full-order state.  However, the errors incurred by this
time-integration procedure (compared with standard time integration of the
reduced-order model) are not assessed or controlled.  This can be problematic
in realistic scenarios, where error estimators and bounds are essential.
Another class of techniques known as \emph{a priori} model reduction methods
\citep{ryckelynck2005phm,skippingSteps} build a reduced-order model `on the
fly', i.e., over the course of a given time integration. These techniques aim
to use the reduced-order model at as many time steps as possible; they revert
to the
high-fidelity model when the reduced-order model is deemed to be inaccurate.
In effect, these techniques employ the reduced-order model as a tool to accelerate the
high-fidelity-model simulation.  In contrast, this work aims to accelerate the
reduced-order-model simulation itself. Further, these methods differ from the
present context in that there are no training experiments conducted \emph{a
priori} from which to glean insight into the model's temporal behavior.

In this work, we propose a method that exploits a set of complete trajectories
observed during training simulations to decrease the temporal complexity of a
reduced-order-model simulation. The method 1) forecasts the unknown variable in
the reduced-order system of nonlinear algebraic equations, and 2) uses this forecast as
an initial guess for the Newton-like solver. To compute the forecast, the
method employs the Gappy POD method \citep{sirovichOrigGappy}, which
extrapolates the unknown variable at future time steps by exploiting the
unknown variable for the previous $\memory$ time steps (where $\memory$ can be
interpreted as the
memory of the process), and a database of time histories of the unknown
variable. If the forecast is accurate, then the Newton-like solver will require very
few iterations to converge, thereby decreasing the number of linear-system
solves needed for the simulation. The method is straightforward to implement: the
(offline) training stage simply requires collecting an additional set of snapshots during
the training simulations. In some scenarios, no additional offline work is
required. The (online) reduced-order-model simulation simply requires an external
routine for determining the initial guess for the Newton-like solver.

% \subsection{Related work}

% \marginpar{ktc - bart and jaideep to do}

%  \begin{itemize} 
%    \item Methods for mitigating time-discretization burden
%    \item Methods for decreasing temporal complexity
%    \item Methods for providing an initial guess to the Newton solver for
%      implicit dynamics
%    \item Methods for forecasting
%    \item Methods for exploiting experimental temporal correlations for
%      surrogate modeling
%      \begin{itemize} 
%         \item \citep{glaz2010reduced}
% 	\item Temporal behavior has been used previously in reduced-order
% 	  modeling. Bos et al.\ \citep{bos2004als} proposed a reduced-order
% 	  model in the context of explicit time integration wherein the
% 	  generalized coordinates are computed based on a best-linear-unbiased
% 	  (BLU) estimate approach. Here, the reduced state coordinates at time
% 	  step $n+1$ are computed using empirically derived correlations between
% 	  the reduced state coordinates and 1) their value at the previous time
% 	  step, 2) the forcing input at the previous time step, and 3) a subset
% 	  of the full-order state.
%      \end{itemize}
%  \end{itemize}

% \subsubsection{Existing forecasting techniques}

\section{Problem formulation}
This section provides the context for this work. Section \ref{sec:fom}
describes the class of full-order models we consider, which includes first-
and second-order ODEs numerically solved by implicit time integration. Section
\ref{sec:rom} describes the reduced-order modeling strategies for which the
proposed technique is applicable.

\subsection{Full-order model}\label{sec:fom}
\subsubsection{First- and second-order ODEs}
First, consider the parameterized nonlinear first-order ODE corresponding to
the full-order model of a dynamical system:
 \begin{gather} \label{eq:ODE1}
 \dot \state = \f \left(\state;t,\forcing\left(t\right),\param\right)\\
 \label{eq:ODE1last}\state(0,\forcing,\param) = \stateInitial.
 \end{gather} 
Here, time is denoted by $t\in\left[0,\lastT\right]$, the
time-dependent forcing inputs are denoted by $\forcing:\left[0,\lastT\right]\rightarrow
\RR{\nforcing}$, the time-independent parametric inputs are denoted by
$\param \in \paramDomain\subseteq\RR{\nparam}$ with $\paramDomain$ denoting
the parameter domain, and $\f :\RR{\ndof }\times \left[0,\lastT\right]
\times \RR{\nforcing}\times \RR{\nparam}\rightarrow \RR{\ndof }$
is nonlinear in at least its first argument. The state is denoted by 
$\state\equiv \state(t,\forcing,\param)\in \RR{\ndof }$ with $\ndof $ denoting the
number of degrees of freedom in the model. The parameterized initial condition
is $\stateInitialNo:\RR{\nforcing}\rightarrow\RR{\ndof }$.

Because this work addresses both first- and second-order ODEs, consider also
the parameterized nonlinear second-order ODE corresponding to the full-order
model of a dynamical system:
 \begin{gather} \label{eq:ODE2}
 \ddot \state = \g \left(\state,\dot\state;t,\forcing\left(t\right),\param\right)\\
 \state(0,\forcing,\param) = \stateInitial\\
 \label{eq:ODE2last}\dot\state(0,\forcing,\param) = \vInit(\param).
 \end{gather} 
Here, the function $\g :\RR{\ndof }\times \RR{\ndof }\times \left[0,\lastT\right]
\times \RR{\nforcing}\times \RR{\nparam}\rightarrow \RR{\ndof }$
is nonlinear in at least its first or second argument, and the parameterized initial
velocity is denoted by $\vInit:\RR{\nforcing}\rightarrow\RR{\ndof }$.\footnote{Note that an $\ndof $-dimensional second-order
ODE can be rewritten as $2N$-dimensional first-order ODE.}

\subsubsection{Implicit time integration}\label{sec:implicitTime}
Given forcing and parametric inputs, the numerical solution to the full-order
model described by Eqs.~\eqref{eq:ODE1}--\eqref{eq:ODE1last} or
\eqref{eq:ODE2}--\eqref{eq:ODE2last} can be computed via numerical integration.
For stiff systems, an implicit integration method is often the
most computationally efficient choice; it is even essential in many cases
\citep{hairer1993solving}. When
an implicit time integrator is employed, $s$ coupled
$\ndof $-dimensional systems of nonlinear algebraic equations are solved at each time step
$n=1,\ldots,\totalTS$, where $\totalTS$ denotes the total number of time
steps:
\begin{equation} \label{eq:fomNonlin}
\Rni\left(\unknownnone,\ldots,\unknownns;\forcing,\param\right)
=0,\quad i = 1,\ldots,s.
\end{equation} 
Here, the function $\Rni:\RR{\ndof }\times\cdots\times
\RR{\ndof }\times\RR{\nforcing}\times\RR{\nparam}\rightarrow \RR{\ndof }$ is
nonlinear in at least one of its first $s$ arguments and the unknowns
$\unknownni\in\RR{\ndof }$, $i=1,\ldots,s$
are implicitly defined by \eqref{eq:fomNonlin}. As discussed in \ref{app:implicitFirst} and
\ref{app:implicitSecond}, the unknowns $\unknownni$ represent the state,
velocity, or acceleration at points $t^{n-1} + c_i\timestep^n$, where
$c_i\in\left[0,1\right]$ is defined by the time integrator:
 \begin{equation} 
\unknownni\equiv \unknownni(\forcing,\param)\equiv
\unknown(t^{n-1}+c_i\timestep^n;\forcing,\param).
  \end{equation} 
Thus, a superscript $n$ denotes the value of a quantity at time $t^n \equiv
\sum\limits_{k=1}^n\timestep^k$, a superscript
$n,i$ denotes the value of a quantity at time $t^{n,i}\equiv
\sum\limits_{k=1}^{n-1}\timestep ^k + c_i\timestep^n$, and $\timestep $ denotes the time-step size.

After the unknowns are computed by solving Eq.~\eqref{eq:fomNonlin}, the
state is explicitly updated as
\begin{equation} \label{eq:stateUp}
\state^n = \gamma \state^{n-1} + \sum_{i=1}^s\delta_i\unknownni,
\end{equation} 
where $\gamma$ and $\delta_i$, $i=1,\ldots,s$ are scalars defined by the
integrator. For second-order ODEs, the velocity is also updated explicitly as
\begin{equation} \label{eq:velUp}
\dot \state^n = \epsilon \dot \state^{n-1} + \sum_{i=1}^s\xi_i\unknownni,
\end{equation} 
where $\epsilon$ and $\xi_i$, $i=1,\ldots,s$ are also scalars defined by the
integrator. 
\ref{app:implicitFirst} and \ref{app:implicitSecond} specify
the form of Eqs.\ \eqref{eq:fomNonlin}--\eqref{eq:velUp} for important classes
of implicit numerical integrators for first- and second-order ODEs,
respectively.  

The chief computational burden of solving Eq.~\eqref{eq:ODE1} with an implicit
integrator lies in solving nonlinear equations \eqref{eq:fomNonlin} at each
time step; this is typically done with a Newton-like method. In particular,
if $\newtonItsAvg$ denotes the average number of Newton-like iterations required to solve
	\eqref{eq:fomNonlin}, then the full-order-model simulation requires solving
	$\newtonItsAvg\totalTS$ linear systems of dimension $sN$.\footnote{Assuming the
	Jacobian of the 
	residual is sparse with an average number of nonzeros per row
	$\sparsity\ll \ndof $, the dominant computational cost of solving Eqs.\
	\eqref{eq:fomNonlin} for the entire simulation is $\mathcal O\left(\sparsity^2sN\totalLinSys\right)$ if a
	direct linear solver is used. It is $\mathcal O\left(\linearIt\sparsity sN \totalLinSys \right)$ if an iterative
	linear solver is used. Here, $\linearIt$ denotes the average number of
	matrix-vector products required to solve each linear system in the case of
	an iterative linear solver.}
	We denote the simulation's \emph{spatial complexity} to be the computational
	cost of solving each linear system; we consider the simulation's
	\emph{temporal complexity} to be the total number of linear-system solves.

The spatial complexity contributes significantly to the computational burden
for large-scale systems because $\ndof $ is large. However, the temporal
	complexity is also significant for such problems. First, the number of total
	time steps $\totalTS$ is often proportional to a fractional power of $\ndof$.
	This occurs because refining the mesh in space often
	necessitates a decrease in the time-step size to balance the spatial and
	temporal errors.\footnote{This is not necessarily true for explicit time-integration schemes, when the time-step size is limited by stability rather than accuracy. In this
case, Krysl et al.\ \citep{krysl_nl_rom_dynam_struct_01} showed that employing
a low-dimensional subspace for the state may improve stability and therefore
permit a larger time-step size.  As a result, the reduced-order state
equations can be solved fewer times than the full-order state equations.}
Second, the average number of Newton-like iterations $\newtonItsAvg$ can be
	large when the problem is highly nonlinear and large time steps are taken,
	which is common for implicit integrators.  Under these conditions, the
	initial guess for the Newton solver, which is often taken to be a polynomial
	extrapolation of the unknown, can be far from the true value of the
	unknown.

In many cases (e.g., linear multi-step methods, single-stage Runge--Kutta
schemes), $s=1$. For this reason, and for the sake of notational clarity, the
remainder of this paper assumes $s=1$, and $\unknown^n$ designates the value of the
unknown variable at time $t^{n,1}$. However, we note that the proposed
technique can be straightforwardly extended to $s>1$.

\subsection{Reduced-order model}\label{sec:rom}

Nonlinear model-reduction techniques aim to generate a low-dimensional model that is
inexpensive to evaluate, yet captures key features of the full-order model.
To do so, these methods first perform analyses of the full-order model for a
set of $\nTrain$ training parametric and forcing points
$\{(\forcing^k,\param^k)\}_{k=1}^{\nTrain}$ during a computationally
intensive `offline' training stage. These analyses may include 
integrating the equations of motion, modal decomposition, etc.

Then, the data generated during these analyses are employed to decrease the
the cost of each linear-system solve
via two approximations: 1) dimensionality reduction, 2) nonlinear-function
approximation (spatial-complexity reduction).  Once these approximations are defined,
the resulting reduced-order model is employed to perform computationally
inexpensive analyses for any inputs
 during the `online' stage.

\subsubsection{Dimensionality reduction}\label{sec:dimRed}
Model-reduction techniques decrease the number of degrees of freedom by
computing an approximate state $\approxstate\approx \state$ that lies in an
affine
%\footnote{Other approaches employ a linear trial subspace
%$\approxstate(t,\forcing,\param) = \podstate\redstate(t,\forcing,\param)$.
%However, this requires projecting the (given) initial condition onto
%$\Range{\podstate}$, which introduces additional error, yet does not decrease
%the dimension or complexity of the model.
%For this
%reason, we employ the affine trial subspace of Eq.~\eqref{eq:redstate}.} 
trial subspace of dimension $\nstate\ll \ndof $:
 \begin{gather} \label{eq:redstate}
 \approxstate(t,\forcing,\param) = \stateRef +
\podstate\redstate(t,\forcing,\param) \\
\label{eq:redvel}\dot\approxstate(t,\forcing,\param) = 
\podstate\dot\redstate(t,\forcing,\param)\\
\label{eq:redacc}\ddot\approxstate(t,\forcing,\param) = 
\podstate\ddot\redstate(t,\forcing,\param).
 \end{gather} 
Here, the trial basis (in matrix form) is denoted by
$\podstate\equiv\vecmat{\podstatevec}{\nstate}\in\RR{\ndof \times \nstate}$ with
$\podstate^T\podstate = \identity$. The
generalized state is denoted by $\redstate\equiv \left[\redstateEntry{1}\ \cdots\
\redstateEntry{\nstate}\right]^T\in\RR{\nstate}$. The reference state is
$\stateRefNo\in\RR{\ndof}$, which is often set to zero. The initial condition for the
reduced-order model is obtained by projecting the full-order-model initial
condition onto this affine subspace such that 
\begin{gather}
\approxstate(0,\forcing,\param)
= \stateRef +
\podstate\podstate^T\left(\stateInitial
- \stateRef\right)\\
\dot\approxstate(0,\forcing,\param)
= \podstate\podstate^T\vInit(\param).
\end{gather}
When the unknown variable computed at each time step (see Section
\ref{sec:implicitTime}) corresponds to the state, velocity, or acceleration,
this dimensionality reduction for the state results in the following
dimensionality
reduction for the unknown:
 \begin{equation} 
\unknownApprox(t,\forcing,\param) = \unknownRef +
\podstate\redUnknownTime(t,\forcing,\param),
  \end{equation} 
	where $\unknownRef= \stateRef$ if the unknown is the state and $\unknownRef=
	0$
	otherwise, and 
	$\redUnknownTime\equiv\left[\redUnknownTimeEntry{1}\ \cdots \
	\redUnknownTimeEntry{\nstate}\right]^T\in\RR{\nstate}$
	denotes the vector of generalized unknowns.

Substituting Eqs.\ \eqref{eq:redstate}--\eqref{eq:redvel}  into
\eqref{eq:ODE1} yields
 \begin{equation} \label{eq:rom1}
 \podstate\dot {\redstate} = \f \left(\stateRef +
 \podstate\redstate;t,\forcing\left(t\right),\param\right),
 \end{equation} 
Alternatively, substituting Eq.\ \eqref{eq:redstate}--\eqref{eq:redacc}  into
\eqref{eq:ODE2} yields
 \begin{equation} \label{eq:rom2}
 \podstate\ddot {\redstate} = \g \left(\stateRef + \podstate\redstate,\podstate\dot\redstate;t,\forcing\left(t\right),\param\right).
 \end{equation} 

The overdetermined ODEs described by \eqref{eq:rom1} and \eqref{eq:rom2} may not be solvable, because
$\text{image}(\f )\not\subset\text{range}(\podstate)$ and
$\text{image}(\g )\not\subset\text{range}(\podstate)$ in general. Several
methods exist to compute an approximate solution.
\paragraph{Project, then discretize in time}
This class of model-reduction methods first carries out a projection process
on the ODE followed by a time-integration of the resulting low-dimensional ODE.
The (Petrov--Galerkin) projection process enforces orthogonality of the
residual corresponding to the overdetermined ODE
\eqref{eq:rom1} or \eqref{eq:rom2} to an $\nstate$-dimensional test subspace
$\text{range}(\testBasis)$, with
$\testBasis\in\RR{\ndof \times \nstate}$. Assuming $\testBasis^T\podstate$ is
invertible, this leads to the following for first-order ODEs:
 \begin{gather} 
 \label{eq:pgODE1}\dot {\redstate} =
 \left(\testBasis^T\podstate\right)^{-1}\testBasis^T \f\left(\stateRef +
\podstate\redstate;t,\forcing\left(t\right),\param\right).
 \end{gather} 
 For second-order ODEs, the result is
 \begin{gather} \label{eq:pgODE2}
 \ddot {\redstate} = \left(\testBasis^T\podstate\right)^{-1}\testBasis^T
 \g\left(\stateRef +
 \podstate\redstate,\podstate\dot\redstate;t,\forcing\left(t\right),\param\right),
 \end{gather} 
Galerkin projection
corresponds to the case where $\testBasis = \podstate$.

Because Eq.\ \eqref{eq:pgODE1} (resp.\ \eqref{eq:pgODE2}) is an ODE of the
same form as \eqref{eq:ODE1} (resp.\ \eqref{eq:ODE2}), it can be solved
using the same numerical integrator that was used to solve Eq.~\eqref{eq:ODE1}
(resp.\ \eqref{eq:ODE2}). Further, the same time-step sizes are often
employed, as the time-step size is determined by accuracy (not stability) for
implicit time integrators. For both first- and second-order ODEs, this again
leads to a system of nonlinear equations to be solved at each time
step
$n=1,\ldots,\totalTS$:
\begin{equation} \label{eq:redStateEq}
 \left(\testBasis^T\podstate\right)^{-1}\testBasis^T\Rn\left(\unknownRef +
\podstate\redUnknownn;\forcing,\param\right)=0.
\end{equation} 
The unknown $\redUnknownn$ can be computed by applying
Newton's method to \eqref{eq:redStateEq}.  Then, the explicit updates
\eqref{eq:stateUp}--\eqref{eq:velUp} can proceed as usual to compute the
resulting state.

\paragraph{Discretize in time, then project}
 This class of model-reduction techniques first applies the same numerical
 integrator that was used to solve \eqref{eq:ODE1} to 
 the overdetermined ODE \eqref{eq:rom1} or \eqref{eq:rom2}. However,
 the resulting algebraic system of $\ndof $ nonlinear equations in $\nstate $
 unknowns remains overdetermined:  
 \begin{equation} \label{eq:discreteOv}
\Rn\left(\unknownRef +
\podstate\redUnknownn
;\forcing,\param\right)=0.
  \end{equation} 
  
To compute a unique solution to Eq.~\eqref{eq:discreteOv}, orthogonality of the
discrete residual $\Rn$ to a test subspace
$\text{range}\left(\testBasis\right)$ can be enforced.  However, this
leads to a reduced system of nonlinear equations equivalent to 
\eqref{eq:redStateEq}. So, in this case, the two
classes of model-reduction techniques are equivalent.

	On the other hand, to compute a unique solution to \eqref{eq:discreteOv},
	the discrete-residual norm can be
	minimized
	\citep{CarlbergGappy,carlbergJCP,bui2008model,bui2008parametric,LeGresleyThesis}, which
	ensures \emph{discrete optimality} \citep{carlbergJCP}:
\begin{equation} \label{eq:resMin}
\redUnknownn=
\arg\min_{\dummy\in\RR{\nstate 
}}\|\Rn\left(\unknownRef +
\podstate\dummy
;\forcing,\param\right)\|^2_2.
\end{equation} 
The unknown $\redUnknownn$ can be computed by applying a
Newton-like nonlinear least-squares method (e.g., Gauss--Newton,
Levenberg--Marquardt) to problem~\eqref{eq:resMin}.  Again, explicit updates for the
state \eqref{eq:stateUp}--\eqref{eq:velUp} can proceed after the unknowns are
computed.

\subsubsection{Spatial-complexity reduction}\label{sec:hyper}

For nonlinear dynamical systems, the dimensionality reduction described in Section
\ref{sec:dimRed} is insufficient to guarantee a reduction in the computational
cost of each linear-system solve. The reason is that the full-order residual
depends on the state, so it must be recomputed and subsequently projected or
minimized at each Newton-like iteration.

For this reason, nonlinear model-reduction techniques employ a
procedure to reduce the spatial-complexity, i.e., decrease the computational
cost of computing and projecting or minimizing the nonlinear
residual. Such techniques are occasionally referred to as
`hyper-reduction' techniques \citep{ryckelynck2005phm}.  In
particular, the class of `function sampling' techniques
replace the full-order nonlinear residual with an approximation $\Rapprox
\approx \residual$ that is inexpensive to compute.  Then, $\Rn\leftarrow
\Rapprox^n$ is employed in \eqref{eq:redStateEq} or \eqref{eq:resMin} to
compute the unknowns $\redUnknownn$.

Methods in this class can be categorized as
follows:
 \begin{enumerate} 
  \item \emph{Collocation approaches}. These methods employ a residual
	approximation that sets many of the residual's entries to zero:
 \begin{equation} \label{eq:collocation}
 \Rapprox^n = \restrict^T\restrict \Rn.
  \end{equation} 
	Here, $\restrict\in\{0,1\}^{\nsample\times \ndof }$ is a sampling matrix
	consisting of $\nsample\ll \ndof $ selected rows of
	$\identity_{\ndof \times \ndof }$. This approach has been developed for Galerkin projection
	\citep{astrid2007mpe,ryckelynck2005phm} and discrete-residual minimization
	\citep{LeGresleyThesis}. 
	\item \emph{Function-reconstruction approaches}. These methods employ a
	residual approximation that computes a few entries of the residual or
	nonlinear function, and subsequently `fills in' the remaining entries via
	interpolation or least-squares regression. That is, these methods apply one
	of the following approximations:
 \begin{align} 
 \label{eq:gappyRes}\tildeRn &=
 \podres\left(\restrict\podres\right)^+\restrict \Rn\\
 \label{eq:gappyFirst}\tilde \f  &= \podf\left(\restrict\podf\right)^+\restrict \f \\
 \tilde \g  &= \podg\left(\restrict\podg\right)^+\restrict \g .
  \end{align} 
	Here, $\podres$,  $\podf$, and $\podg$ are empirically derived bases used to
	approximate the nonlinear residual, velocity, and acceleration,
	respectively. A
	superscript $+$ denotes the Moore--Penrose pseudoinverse. When the bases are
	computed via POD, this technique is known as Gappy POD
	\citep{sirovichOrigGappy}. This approach has been developed for Galerkin projection
	\citep{astrid2007mpe,bos2004als,chaturantabut2010journal,galbally2009non,drohmannEOI}
	and discrete-residual minimization \citep{CarlbergGappy,carlbergJCP}. In
	particular, the discrete empirical interpolation (DEIM)
	method \citep{chaturantabut2010journal} is a specific case of Gappy POD for
	first-order ODEs, Galerkin projection, and the interpolatory case, i.e.,
	DEIM uses approximation \eqref{eq:gappyFirst} in Eq.~\eqref{eq:pgODE1} with
	$\testBasis = \podstate$ and
	sets the number of sample indices $\nsample$ equal to the number of columns
	in the matrix $\podf$. The GNAT method \citep{CarlbergGappy,carlbergJCP}
	employs Gappy POD of the residual in a discrete residual minimization
	setting, i.e., GNAT uses approximation \eqref{eq:gappyRes} in
	Eq.~\eqref{eq:resMin}.

	 \end{enumerate}

\section{Temporal-complexity reduction}
While the model-reduction approaches described in the previous section decrease the
computational cost of each linear-system solve (i.e., spatial complexity),
they do not necessarily decrease the number of linear-system solves (i.e., temporal
complexity). The goal of this work is devise a method that decreases this
temporal complexity while introducing no additional error. 

\subsection{Method overview}\label{sec:methodOveriew}

The main idea of the proposed approach is to compute an accurate forecast of
the generalized unknowns at future time steps using the Gappy POD procedure,
and employ this forecast as an initial guess for the Newton-like solver at
future time steps.  

Gappy POD is a technique to reconstruct vector-valued data that has `gaps,'
i.e., entries with unknown or uncomputed values.  Mathematically, the approach
is equivalent to least-squares regression in one discrete-valued variable
using empirically computed basis functions. It was introduced by Everson and
Sirovich~\citep{sirovichOrigGappy} for the purpose of image reconstruction.
It has also been used for static \citep{bui2003proper,bui2004aerodynamic} and
time-dependent \citep{venturi2004gappy,willcox2006ufs} flow field
reconstruction, inverse design~\citep{bui2004aerodynamic}, design variable
mapping for multi-fidelity optimization~\citep{robinson:smo}, and for
decreasing the spatial complexity in nonlinear model
reduction~\citep{astrid2007mpe,bos2004als,CarlbergGappy,carlbergJCP}.  This
work proposes a novel application of Gappy POD: as a method for forecasting
the generalized unknown at future time steps during a reduced-order-model
simulation.

\begin{figure}[htbp] \centering
\subfigure[Offline: the computed time-evolution POD basis for a generalized unknown.]
{\includegraphics[width=.5\textwidth]{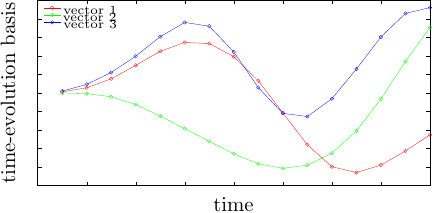}\label{fig:timeBasis}}
\subfigure[Online: time steps taken so far (red), recent time steps used
to compute
forecast (green), forecast (blue)]
{\includegraphics[width=.5\textwidth]{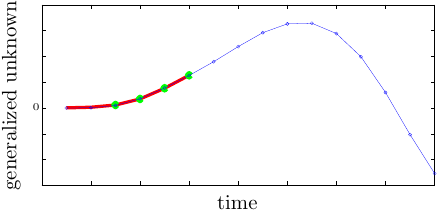}\label{fig:forecast}}
\caption{Graphical depiction of the proposed method}
\end{figure}

During the offline stage, the proposed method computes a `time-evolution basis' for
each generalized unknown $\redUnknownTimeEntry{\basisIndex}$,
$\basisIndex=1,\ldots,\nstate$. Each basis represents the complete time-evolution of a
generalized unknown as observed during training simulations. Figure
\ref{fig:timeBasis} depicts this idea graphically, and Section
\ref{sec:timeBasis} describes a computationally inexpensive way to compute
these bases.

During the online stage, the method computes a forecast of the generalized
unknowns at future time steps via Gappy POD. This forecast employs 1) the
time-evolution bases and 2) the generalized unknowns computed at several previous time
steps. Figure \ref{fig:forecast} depicts this, and Section \ref{sec:forecast}
describes the forecasting method in detail. At future time steps, this forecast is
employed as an initial guess for the Newton-like solver. If the forecast is
accurate, the Newton-like solver will converge in very few iterations; if it
is inaccurate, the Newton-like solver will require more iterations for
convergence. Note that the accuracy of the solution is not hampered in either
case (assuming a globalization strategy is employed). If the number of Newton
iterations required for convergence is large, this indicates an inaccurate initial
guess. When this occurs, the method computes a new forecast using the most
recently computed generalized unknowns.

%\begin{figure}[htbp] \centering
%\includegraphics[width=.55\textwidth]{figureGeneration/timeSnapsOnline3state_train_tr16_sz11_ts33_dt300_sn2_PODcoords}
%\caption{Online: time steps taken so far (red), recent time steps used
%to compute
%forecast (green), forecast (blue)}\label{fig:forecast}
%\end{figure}

The proposed method is expected to be effective if 1) the temporal behavior of
the generalized unknowns is similar across input variation and 2) the original
model is not too weakly nonlinear at each time step. The latter issue can
hamper the proposed method's performance because it is difficult to reduce the
number of Newton iterations if the original number is already very small.
This situation can occur, for example, if the simulation employs a very small time step.
However, this is uncommon for (unconditionally stable) implicit time
integrators, where taking the largest time step while maintaining accuracy is
typically the most computationally efficient approach.

The proposed method is
independent of the dimensionality-reduction or spatial-complexity-reduction scheme employed by
the reduced-order model; further, the method is applicable (without
modification) to both first- and second-order ODEs. The next sections
describe the offline and online steps of the methodology in detail.

\subsection{Offline stage: compute the time-evolution bases}\label{sec:timeBasis}

The objective of the offline stage is to compute the
time-evolution bases that will be used for the online forecast.
Ideally, the bases should be able to describe the time evolution of the
generalized state for any forcing inputs $\forcing$ and parametric inputs
$\param$. If the bases are `bad', then the forecasting step of the algorithm
will be inaccurate, and there may be no reduction in the average number of
Newton-like iterations.

We propose employing a POD basis for the time evolution of the generalized
unknown. This basis is computed \emph{a priori} during `offline' simulations
of the reduced-order model in three steps:
 \begin{enumerate} 
 \item\label{step:unknownSnaps} Collect snapshots of the unknown during each of the $\nTrain$ training
 simulations:
 \begin{equation} 
\snapsMat{k} = \left[\unknown^{0}\left(\forcing^k,\param^k\right)\
\cdots\
\unknown^{\totalTS-1}\left(\forcing^k,\param^k\right)\right]
 \end{equation} 
 for $k=1,\ldots,\nTrain$, with $\snapsMat{k}\in\RR{\ndof \times
	 \totalTS}$.
Here, $\forcing^k\in \RR{\nforcing}$ denotes the forcing inputs for
training simulation $k$,
and $\param^k\in\RR{\nparam}$ denotes the parametric inputs for
training simulation $k$.

 \item \label{step:genUnknownSnaps}Compute the corresponding snapshots of the generalized unknown:
 \begin{align} 
\redSnapsMat{k}&\equiv\podstate^T \left[\snapsMat{k} - \unknownTimeInitk
\ones^T\right] \\
&= \left[
\redUnknownTime^{0}\left(\forcing^k,\param^k\right)\ \cdots\
\redUnknownTime^{\totalTS-1}\left(\forcing^k,\param^k\right)
\right]
 \end{align} 
 for $k=1,\ldots,\nTrain$, where orthogonality of the trial basis
	 $\podstate^T\podstate = \identity$ has been used. Here, $\redSnapsMat{k}\in\RR{\nstate\times
	 \totalTS}$ and $\ones\in\RR{\totalTS}$ denotes a vector of ones.
	 \item\label{step:genUnknownSVD} Compute the time-evolution bases via the
	 (thin) singular value
	 decomposition (SVD). Defining the $j$th column of
	 $\redSnapsMat{k}^T$ as $\redSnapsVec{k}{j}\in\RR{\totalTS}$, $j=1,\ldots,\nstate$, we note
	 that 
	 $\redSnapsVec{k}{j}$ can be interpreted as a snapshot of
	 the time evolution of the $j$th generalized unknown $\redUnknown_j$ during training
	 simulation $k$. Then, this step amounts to
 \begin{gather} 
 \left[\redSnapsVec{1}{j}\ \cdots\ \redSnapsVec{\nTrain}{j}\right]
=\boldsymbol U_j\boldsymbol \Sigma_j \boldsymbol V_j^T\\
\label{eq:timePod}\timebasisj = \vecmatsub{\boldsymbol u}{j}{\ntimej},
 \end{gather} 
 for $j=1,\ldots,\nstate$. Here,
	 $\boldsymbol U_j\equiv\vecmatsub{\boldsymbol u}{j}{\nTrain}\in\RR{\totalTS\times \nTrain}$ and $\ntimej \leq
	 \nTrain$.
 \end{enumerate}
After the time-evolution bases $\timebasisj\in\RR{\totalTS\times \ntimej}$,
$j=1,\ldots,\nstate$ have been computed during the offline stage, they can
be used to accelerate online computations via forecasting. The next section
describes this.

\begin{remark}
In some cases, many of the above offline steps are already completed as part
of the existing model-reduction process. For example, the snapshot matrices
$\snapsMat{k}$, $k=1,\ldots,\nTrain$ in Step \ref{step:unknownSnaps} are
already available if proper orthogonal decomposition (POD) is employed to
compute $\podstate$ and the  time integrator's unknown is the state 
(e.g., linear multistep schemes).
If additionally $\nTrain = 1$ and the POD basis is computed
via the SVD of the reference-centered state
snapshots, i.e.,
$\left[\state^{0}\left(\forcing^1,\param^1\right) - \stateRefArg{\param^1}\
\cdots\
\state^{\totalTS-1}\left(\forcing^1,\param^1\right)- \stateRefArg{\param^1}\right]=\bar{\boldsymbol
U}\bar{\boldsymbol \Sigma} \bar{\boldsymbol V}^T$ with $\podstatevec_i
=\bar{\boldsymbol u}_i$, $i=1,\ldots,\nstate$, then $\redSnapsMat{1}$ of Step
\ref{step:genUnknownSnaps} is
already available as $\redSnapsMat{1} =
\bar{\boldsymbol \Sigma}[1:\nstate,1:\nstate] \bar{\boldsymbol
V}[1:\totalTS,1:\nstate]^T$. Here,  the square bracket indicates a
submatrix over the specified range of row and column indices and $\bar{\boldsymbol
U}\equiv\vecmat{\bar{\boldsymbol u}}{\totalTS}$. Further, in this
case the
matrices $\boldsymbol U_j$ in Step \ref{step:genUnknownSVD} are also 
available
as  $\boldsymbol U_j =\boldsymbol u_{j,1} = \bar{\boldsymbol v}_j$,
$j=1,\ldots,\nstate$, where $\bar{\boldsymbol
V}\equiv\vecmat{\bar{\boldsymbol v}}{\totalTS}$.
\end{remark}
\subsection{Online stage: forecast}\label{sec:forecast}

During the online stage, the method employs a forecasting procedure to define
the initial guess for the Newton-like solver. To compute this forecast, it
uses the time evolution bases (computed offline), and the values of the
generalized unknown at the previous $\memory$ time steps (computed online).
Here, $\alpha$ is considered the `memory' of the process. Because the forecast
is defined at all time steps (see the blue curve in Figure
\ref{fig:forecast}), it is used as the initial guess at future time steps
until the number of Newton iterations exceeds a threshold value
$\newtonThreshold$. This indicates a poor forecast. In this case, the forecast
is recomputed using the most recent values of the generalized unknown. 

If the forecast is accurate, then the number of iterations needed to converge
from the (improved) initial guess will be drastically reduced, thereby
decreasing $\newtonItsAvg$ and hence the temporal complexity.  Algorithm
\ref{alg:forecast} outlines the proposed technique. 

\begin{algorithm}[btp]
\caption{Online: implicit time integration with the temporal-complexity-reduction method}
\begin{algorithmic}[1]\label{alg:forecast}
\REQUIRE Time-evolution bases $\timebasisj\in\RR{\totalTSun\times \ntimej}$,
$\basisIndex = 1,\ldots,\nstate$; maximum memory
$\memoryMax$ with \\ $\memoryMax \geq \max\limits_j\ntimej$; Newton-step threshold $\newtonThreshold$
\ENSURE Generalized state at all $\totalTS$ time steps: $\redstate^n$,
$n=1,\ldots,\totalTS$.\\ Generalized velocity at all $\totalTS$ time steps if
solving a second-order ODE: $\dot\redstate^n$, $n=1,\ldots,\totalTS$.
\FOR[time-step loop]{$n=1,\ldots,\totalTS$}\label{step:forloop}
\IF{forecast $\forecastTime(t^{n-1}+c_1\timestep^n)$ is available}
\STATE Set initial guess for Newton solver to
$\redUnknownTimeZeroEntry{j} = \forecastTimej(t^{n-1} + c_1\timestep^n)$, $j =
1,\ldots,\nstate$.
	\ELSE
	\STATE Use typical initial guess for Newton solver (e.g., polynomial
	extrapolation of unknown).
\ENDIF
\STATE Compute generalized unknowns $\redUnknownn$ by solving reduced-order equations \eqref{eq:redStateEq} or
\eqref{eq:resMin} with a Newton-like method and specified initial guess
$\redUnknown^{n(0)}$.\\
 Let $\newtonIts^n$ denote the number of
Newton-like iterations required for convergence at time step $n$.
\STATE Compute the generalized state $\redstate^n$ using
explicit update \eqref{eq:stateUp}. If solving a second-order ODE, also update the
generalized velocity $\dot\redstate^n$ using explicit update \eqref{eq:velUp}.
\IF[recompute forecast using most recent data]{$\newtonIts^n>\newtonThreshold$ and $(n-1) \geq
\max\limits_j\ntimej$} 
\STATE Set memory $\memory \leftarrow \min(n-1,\memoryMax)$.
\STATE\label{step:foreFirst}\label{step:computeZ} Compute forecasting coefficients $\redtimej$, $j=1,\ldots,\nstate$ using the
unknown at the previous $\memory$ time steps by solving Eq.\
\eqref{eq:gappyTime}.
\STATE\label{step:foreLast} Set forecast to be $\underline{\mathbf{\hat w}}_\basisIndex
=\timebasisj \redtimej$ and define $\forecastTimej\equiv
\unrollinv{\underline{\mathbf{\hat w}}_\basisIndex}$, $j=1,\ldots,\nstate$.
\ENDIF
\ENDFOR
\end{algorithmic}
\end{algorithm}

To compute the forecasting coefficients in step \ref{step:computeZ} of
Algorithm \ref{alg:forecast}, we propose using the Gappy POD
approach introduced by Everson and Sirovich~\citep{sirovichOrigGappy}. This approach
computes coefficients $\redtimej$ via the following linear least-squares
problem:
 \begin{equation} \label{eq:gappyTime}
 \redtimej = \arg\min_{\redtime\in\RR{\ntimej}}\|\restrictn \timebasisj \redtime - \restrictn
 \unroll{\redUnknown_j}\|
 \end{equation} 
Here, the matrix
$\restrictn\in\{0,1\}^{\alpha \times \totalTSun}$ is the sampling matrix that
selects entries corresponding to the previous $\memory$ time steps:
 \begin{equation} 
 \restrictn \equiv \left[\boldsymbol e_{n-\alpha-1}\ \cdots\ \boldsymbol e_{n-1} \right]^T,
  \end{equation} 
where $\boldsymbol e_i$ denotes the $i$th canonical unit vector. Note that $\memory \geq
\ntimej$ is required for Eq.~\eqref{eq:gappyTime} to have a unique solution.
The function $\unrollNone$ in \eqref{eq:gappyTime} `unrolls' time
according to the time discretization; we define $\unrollNone:x\mapsto \mathbf
x$ with $\mathbf x\equiv \left[\text{x}_1\ \cdots
\text{x}_{\totalTSun}\right]^T\in\RR{\totalTSun}$ as
 \begin{equation} 
\text{x}_{n} = x(t^{n-1} + c_1\timestep^n),\quad n=1,\ldots,
M.
  \end{equation} 

The online cost to compute this forecast is very small, as it
entails solving $\nstate$ small-scale linear least-squares problem \eqref{eq:gappyTime} characterized by a
$\memory \times \ntimej$ matrix. For this reason, it is generally advantageous
to employ a small value of $\newtonThreshold$ (i.e., 0 or 1), which results in
a frequent (inexpensive) recomputation of the forecast.

\section{Numerical experiments}

These numerical experiments assess the performance of the proposed
temporal-complexity-reduction method on a structural-dynamics example using
three reduced-order models: Galerkin projection (Eq.\ \eqref{eq:pgODE2} with
$\testBasis = \podstate$), Galerkin projection with least-squares
reconstruction of the residual (Eq.\ \eqref{eq:gappyRes}), and a
structure-preserving reduced-order model \citep{carlberg2012spd}.  We do not
present results for a collocation ROM (see Eq.~\eqref{eq:collocation}), as this approach was
unstable in most cases, even when 60\% of the degrees of freedom were selected
as sample indices (i.e., $\nsample/\ndof = 0.6$).  Section \ref{sec:problemDesc} provides a description of
the problem---a parameterized, damped clamped--free truss structure
subjected to external forces---and details the experimental setup.  We then
consider a sequence of problems that poses increasing difficulty to the method. 

Section \ref{sec:limiting} considers the ideal scenario for the method: the
online points are identical to the training points, and the reduced bases are
not truncated. In this case, the temporal behavior of the system is perfectly
predictable, because (in exact arithmetic) the online response is the same as
the training response. Therefore, we expect the proposed method to work
extremely well.

Section \ref{sec:FREElow2} assesses the method's performance in a more
challenging setting. Here, the online points differ from the training
points (i.e., a predictive scenario), so the temporal behavior is not identical
to that observed during the training simulations. The parametric inputs
correspond to shape parameters and the initial displacement. The external force
is set to zero, which leads to a damped free-vibration problem.  As a result, the
dynamics encountered in this example are relatively smooth.

Section \ref{sec:FRCSTRlow2} considers a more challenging predictive scenario
wherein rich dynamics---generated from a high-frequency external
force---characterize the response. Here, additional parametric inputs are
considered, which correspond to
the magnitudes and frequencies of the high-frequency forces.

Section \ref{sec:FRCSTR2} increases the predictive difficulty, as
the allowable range of the parametric inputs is doubled, leading to a more significant
variation in the responses.

Finally, Section \ref{sec:avg} summarizes the proposed forecasting method's
performance over all experiments and tested reduced-order models.

\subsection{Problem description}\label{sec:problemDesc}

\begin{figure}[htbp] \centering
\includegraphics[width=.5\textwidth]{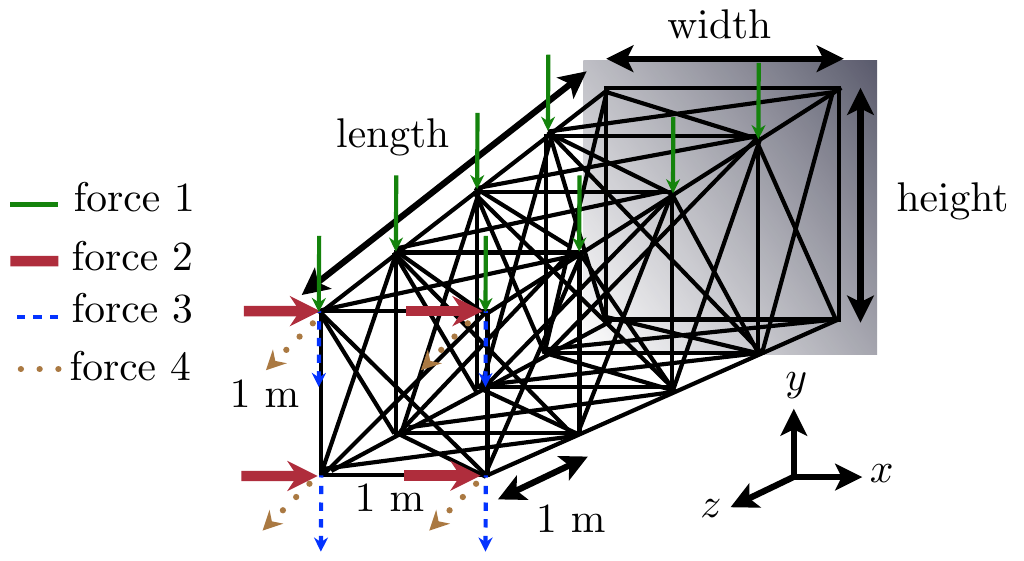}
\caption{Clamped--free parameterized truss structure}
\label{fig:structure}
\end{figure}

Figure~\ref{fig:structure} depicts the parameterized, non-conservative
clamped--free truss structure we consider.
The truss is parameterized by $\nparam = 16$ parametric
inputs $\param\equiv\left(\parami{1},\ldots,\parami{16}\right)\in\paramDomain
= \left[-0.5,0.5\right]^{16}$
that affect the geometry, initial condition, and applied force as described in
Table \ref{tab:parameterization}.
 \begin{table}[htd] 
 \centering 
 \scriptsize{\begin{tabular}{|c|c|c|c|c|c|c|} 
  \hline 
\multirow{3}{*}{length (m)} & bar & \multirow{3}{*}{width (m)}& \multirow{3}{*}{height (m)}&initial
condition
& external-force & external-force  \\
& cross-sectional & & & max magnitude (N)
&magnitude & frequency
 \\
&area ($\mathrm{m}^2$) & & &$\initialCondMag_i$, $i=1,\ldots,4$ & $\forceMag_i$, $i=1,\ldots,4$ &$\lambda_i$, $i=1,\ldots,4$ \\
\hline
$200 + 50\parami{1}$ & $0.0025(1 + 0.5\parami{2})$ & $10(1 + \parami{3})$ & $10(1 +
\parami{4})$&$\fNomi(1
+ 0.5 \parami{i+4})$ & $\fNomi (1 + 0.5\parami{i+8})$& $3\omega_0(1 +
0.5\parami{i+12})$ \\
\hline 
  \end{tabular} }
	\caption{Effect of parameters on truss geometry, initial conditions, and
	applied forces. Here, $\fNomi$, $i=1,\ldots,4$ denote the nominal force
	magnitudes (specified within each experiment) and $\omega_0$ denotes the
	lowest-magnitude eigenvalue at the nominal point
	$\paramNom$.}
\label{tab:parameterization}
	\end{table}
We set the material properties to those of aluminum, i.e., density
$\rho = 2700\ \mathrm{kg/m}^3$ and 
elastic modulus $E=62\times 10^9$ Pa.
The external force is composed of four components:
 \begin{equation} \label{eq:externalForceOne}
\fext(\param,t) = \sum_{i=1}^4\forceTime_i(\param,t)\forceDist_i,
 \end{equation} 
where 
$\forceDist_i\in\RR{\ndof}$, $i=1,\ldots, 4$ correspond to unit loads uniformly distributed
across designated nodes and $\forceTime_i:\paramDomain\times
\left[0,\lastT\right]\rightarrow\RR{}$, $i=1,\ldots, 4$
denote the $\nforcing = 4$ forcing inputs. Figure~\ref{fig:structure} depicts the spatial distribution of the forces, which lead
to vectors $\forceDist_i$, $i=1,\ldots, 4$ through the finite-element
formulation described below. The parameterized, time-dependent magnitudes of these forces
are
\begin{gather} \label{eq:externalForceTwo}
\forceTime_i(\param,t) = 
\begin{cases}
\forceMag_i\left(\param\right)\sin\left(\lambda_i(\param)\left(t - \lastT/4\right)\right),\quad
t\geq \lastT/4\\
0,\quad \mbox{otherwise}
\end{cases}
,
\end{gather} 
where $\forceMag_i:\paramDomain\rightarrow \RR{}$ and
$\lambda_i:\paramDomain\rightarrow \RR{}$, $i=1,\ldots, 4$ denote the maximum force
magnitudes and force frequencies, respectively.
Similarly, the initial displacement is composed of
four components
 $
\stateInitial = \sum_{i=1}^4\initialCondMag_i(\param)\initialCond_i,
 $
 where $\initialCond_i$ is the steady-state displacement of the truss subjected
 to load $\forceDist_i\forceMag_i\left(\paramNom\right)$ with $\paramNom =
 \left(0,\ldots,0\right)$ denoting the nominal point in parameter space. 
 The initial velocity is set to zero $\vInit = 0$, and
 the reference
 configuration is simply the undeformed truss (in equilibrium) represented by
 $\stateRefNo = 0$.

% fixed
%net_force_top = inputs(18) * 1000 * 9.81; %1000 kg about 1 ton
	% inputs(18) is in the name of the file
%net_force_end_z = net_force_top
%net_force_end = 0.2  * net_force_top
% max_frequency = smallest_eigenvalue * 3.0
% magnitude = max_magnitude * sin(freq * (time - delay));
% delay  = 0.25

% INITIAL CONDITION: for NOMINAL configuration
% 1) D = (0.5 * MU + 1) * DISPLAC
% 2) D = (0.5 * MU + 1) * DISPLAC
% 3) D = (0.5 * MU + 1) * DISPLAC
% 4) D = (0.5 * MU + 1) * DISPLAC

% FORCING
%   5) max_magnitude = [(force_per_node_length) * (1 + 0.5 * inputs(8)), ...
%   6)  (force_per_node_end) * (1 + 0.5 * inputs(9)), ...
%   7)  (force_per_node_end) * (1 + 0.5 * inputs(10)),...
%   8)  (force_per_node_end_z) * (1 + 0.5 * inputs(11)),...
%   ];
% 9)--12) frequency = (0.5 * inputs(12:15) + ones(size(inputs(12:15)))) * max_frequency;

%GEOMETRY
%MU2 13) Surface area: 0.0025 m^2 (1 + 0.5 * inputs4)
%MU3 14) a = 10m * (1 + 0.5 * inputs(5)); % width at base
%MU4 15)c = 10m * (1 + 0.5 * inputs(6)); % height at base
%MU1 16) length = 200m + 50 * mu_16

The problem is discretized by the finite-element method. The model consists
of 16 three-dimensional bar elements per bay with three degrees of
freedom per node; this results in 12 degrees of freedom per bay. We consider a
problem with 250 bays, and therefore $\ndof=3\times 10^3$ degrees of freedom in the
full-order model. The bar elements model geometric nonlinearity, which results
in a high-order nonlinearity in the internal force. 
This discretization results in the following equations of motion for the
full-order model:
 \begin{equation} \label{eq:sdmEom}
 \Mparam\ddot \state + \Cparam\dot \state + \fintparam = \fext(t;\param).
  \end{equation} 
Here, $\Mparam\in\RR{\ndof\times \ndof}$ denotes the symmetric-positive-definite mass
matrix, the internal force is denoted by
$\fint:\RR{\ndof}\times\paramDomain\rightarrow\RR{\ndof}$, and
the symmetric-positive-semidefinite Rayleigh viscous damping matrix,
denoted by $\Cparam\in\RR{\ndof\times \ndof}$, is of the form
 \begin{equation}\label{eq:dissFun}
 \Cparam = \rayleighMass\Mparam + \rayleighStiff\nabla_\state
 \fintparamArg{\stateInitialNo}.
  \end{equation} 
Note that $\nabla_\state \fintparamArg{\stateInitialNo}$ represents the
tangent stiffness matrix at the initial condition.
Here, $\alpha$ and
$\beta$ are chosen such that the damping ratio is $\zeta = 15\deg$
%$\sin(5 \pi/180)$ 
for the uncoupled ODEs associated with the
smallest two eigenvalues of the matrix pencil
$\left(\M(\paramNom),\nabla_{\state}\fintparamArgs{0}{\paramNom}\right)$
 \citep{chowdhury2003computation}.

The equations of motion \eqref{eq:sdmEom} can be rewritten in the standard
form of Eqs.\ \eqref{eq:ODE2}--\eqref{eq:ODE2last} as
 \begin{gather} \label{eq:sdmEomStandard}
 \ddot\state = \Mparam^{-1}\left(\fext(t;\param) -\Cparam\dot \state -
 \fintparam\right)\\
 \state(0,\forcing,\param) = \stateInitial\\
 \label{eq:sdmEomStandardLast}\dot\state(0,\forcing,\param) = \vInit(\param).
  \end{gather} 
	The nonlinear function defining the acceleration for the second-order ODE is then
 \begin{equation} 
	\g\left(\state,\dot\state;t,\forcing,\param\right) =
\Mparam^{-1}\left(\fext(t;\param) -\Cparam\dot \state -
 \fintparam\right).
  \end{equation} 
	 
	We employ an implicit Nystr\"om time integrator with constant timestep size
	$\timestep = \timestep^n$, $n=1,\ldots, \totalTS$ to compute the numerical
	solution to Eqs.\ \eqref{eq:sdmEomStandard}--\eqref{eq:sdmEomStandardLast} in the time interval $[0,\lastT]$ with $\lastT=25$ seconds.  In particular, we employ the implicit midpoint rule for both partitions.
	This leads to discrete equations \eqref{eq:nystrom} to be solved at each
	time step with explicit updates
	\eqref{eq:nystromState}--\eqref{eq:nystromVel} and parameters $s=1$, $\hat
	a_{11} = 1/2$, $\bar a_{11} = 1/4$, $\hat b_1 = 1$, $\bar b_1 = 1/2$, $c_1 =
	1/2$. The unknowns are equivalent to the acceleration at the half time
	steps: $\unknown^n =\ddot\state\left(t^{n-1}+1/2h\right)$, $n=1,\ldots,\nt$.
	Multiplying the corresponding residual by $\Mparam$ yields
\begin{equation} \label{eq:midpoint}
\Rn\left(\unknown^n\right) =  \Mparam\unknown^n +
\Cparam\left[\dot \state^{n-1} + \frac{1}{2}h
\unknown^n\right] +
 \fintparamArg{\state^{n-1} + \frac{1}{2}h\dot \state^{n-1}+\frac{1}{4}h^2 
\unknown^n}-
\fext(t^{n-1} +\frac{1}{2} h;\param).
\end{equation} 
To solve $\Rn\left(\unknown^n\right) = 0$ at each time step, 
We employ a globalized Newton solver with a More--Thuente
linesearch \citep{poblano}. Except when noted, convergence of the Newton iterations is declared
when the residual norm reaches $10^{-6}$ of its value computed using a zero
acceleration and the values of the displacement and velocity at the beginning
of the timestep.  The linear system arising at each Newton iteration is solved
directly.

The experiments compare the performance of three reduced-order models: Galerkin
projection (Eq.~\eqref{eq:redStateEq} with $\testBasis = \podstate$), Galerkin
projection with Gappy POD residual
approximation (Eq.~\eqref{eq:gappyRes}), and
a model-reduction method based that preserves the classical Lagrangian structure
intrinsic to the problem (Ref.~\citep{carlberg2012spd}, proposal 1).  To construct the
reduced-order models, we collect snapshots of the required quantities for
$\param\in \Dtrain\subset \paramDomain$ and $t\in\left[0,\lastT\right]$.  The trial
basis $\podstate$ is determined via POD. We collect snapshots of the state
\begin{equation} \snaps{\state}=\{\state^{n-1} + h\dot \state^{n-1} +
\frac{h}{2}\ddot \state^{n,1} \ | \ n=1,\ldots,\nt;\ \param\in\Dtrain\}
\end{equation} and set the trial basis to
$\podstate=\podArgs{\snaps{\state}}{\nu_\state}$, where
$\nu_\state\in\left[0,1\right]$ is an `energy criterion' and $\podArgsNo$ is
defined by Algorithm \ref{PODSVD} in \ref{app:POD}. The reference state is set
to $\stateRefNo = 0$, as this is the equilibrium state for this problem 
\citep{carlberg2012spd}. For Galerkin projection
with least-squares (Gappy POD) residual reconstruction, the following snapshots are
collected during the (full-order model) training simulations: 
\begin{equation}
\snaps{\residual}=\{\Rn\left(\unknown^{n(k)}\right)\ |\ n=1,\ldots,\nt;\
k=0,\ldots,K^n-1;\ \param\in\Dtrain\}.  \end{equation} 
Here, $K^n$
denotes the number of Newton steps taken at time step $n$. The residual basis 
is set to $\podres=\podArgs{\snaps{\residual}}{\nu_{\residual}}$ with
$\nu_{\residual}\in[0,1]$.
For the structure-preserving method, we also collect snapshots of both the mass matrix and the
external forcing vector:
\begin{gather}
\snaps{\Mass}=\{\Mparam\ |\ \param\in\Dtrain\}\\
\snaps{\fext}=\{\fext(t^n;\param)\ |\ n=1,\ldots,\nt;\
\param\in\Dtrain\}.
\end{gather} 
The POD basis for the external force employed by the structure-preserving
method is set to $\podstate_{\fext}=\podArgs{\snaps{\fext}}{\nu_{\fext}}$ with
$\nu_{\fext}\in[0,1]$

%The same
%sampling matrix $\restrict$ is used for the Gappy POD and structure-preserving
%approximations; it is determined using  the GNAT model-reduction method's
%approach for selecting the sample matrix \citep{CarlbergGappy,carlbergJCP}.

Reduced-order models with spatial-complexity reduction employ the same sampling matrix
$\sampleMat$, which is generated using GNAT's greedy sample-mesh algorithm
\citep[Algorithm 3]{carlbergJCP}.\footnote{Greedy-algorithm parameters are
$\podres = \podjac =
\podArgsNo_{\mathbf{r}}$ a POD basis computed using Algorithm \ref{PODSVD} with
snapshots of the numerical residual over
all timesteps and Newton iterations during the full-order-model training
simulations and an energy criterion of
$
\energyCrit\leftarrow \energyCrit_\mathbf{r} = 1-10^{-2} 
$, a target number of sample nodes $n_s =
\nsample/\nu$ with $\nu = 3$ unknowns per node (the $x$-,
$y$-, and $z$-displacements),
an empty seeded sample-node set $\mathcal N =
\emptyset$, and $n_c$ equal to the number of columns in $\podArgsNo_{\mathbf{r}}$.} These models are also implemented using the sample-mesh
concept \citep[Section 5]{carlbergJCP}.  
For the structure-preserving method \citep{carlberg2012spd}, we solve the reduced-basis-sparsification
unconstrained optimization problem using the Poblano toolbox
\citep{poblano}.\footnote{The
initial guess for each of these problems is chosen as
$\restrict^T\restrict \podstate$.}

In all experiments, the proposed forecasting method employs untruncated
time-evolution bases: $\ntimej=\nTrain$, $j=1,\ldots, \nstate$. We compare
its performance with that of the most common approach for generating an
initial guess: a polynomial
extrapolation of varying degree. Note that polynomial
extrapolations of different degrees employ a different number of previous
solutions to generate an initial guess; in our experiments, we associate the
number of previous solutions employed with a `memory' $\memory$. For example,
a zeroth-order extrapolation requires the previous solution
($\unknown^{n(0)} = \unknown^{n-1}$), so $\memory = 1$ in this case. When no
previous solution is used (i.e.,
$\memory = 0$), the polynomial-extrapolation approach uses $\ddot
\state^{n,1} = \unknown^{n(0)} = 0$. In all experiments, the full-order model
employs a zeroth-order extrapolation for the initial guess.

The output of interest is the $y$-displacement of the bottom-left node of the
end face of the truss in Figure \ref{fig:structure}.  We denote this
(parameterized, time-dependent) quantity by $d\in\RR{}$.
To quantify
the performance of the reduced-order models, the following
metrics are used:
\begin{gather}
\efomone=\frac{\frac{1}{\nt}\sum\limits_{n=0}^{\nt}
|d^n-d^n_{\FOM}|}{\max\limits_nd^n_\FOM-\min\limits_nd^n_\FOM}\\
\temporalRedFactor =\frac{\newtonItsAvg_\FOM}{\newtonItsAvg} \\
\speedup = \frac{\walltime_\FOM}{\walltime}
\end{gather}
Here, error measure $\efomone$ designates the scaled $\ell_1$ norm of the discrepancy in the
output predicted by a reduced-order model. The temporal-complexity-reduction
factor is denoted by  $\temporalRedFactor$, where $\newtonItsAvg$ denotes the
average number of Newton-like steps taken per time step over the course of a
simulation.  The
speedup is denoted by $\speedup$ with $\walltime$ denoting the wall time
required for a simulation.  A subscript `FOM' denotes a quantity computed
using the full-order model.

All computations are carried out in
Matlab on a Mac Pro with 2 $\times$ 2.93 GHz 6-Core Intel Xeon processors and
64 GB of memory.
%snapshots = -1
	% snapshot 0
	% project initial condition

\subsection{Ideal case: unforced, invariant inputs, no truncation of
bases}\label{sec:limiting}
%\subsection{FREE2low2 nominal and ideal}
	%can use larger time step!

This experiment	explores the ideal case for the method: the
online inputs equal the training inputs, and the bases are not truncated
($\nu_\state = \nu_{\residual} = 1.0$). The resulting basis dimensions are $\nstate
= 100$ for the reduced basis and $329$ for the residual basis (i.e.,
$\podf\in\RR{\ndof\times 329}$). In this scenario, the full-order model's
temporal behavior encountered online is exactly the same as that observed
during training simulation; for this reason, we expect the proposed method to
perform very well. We consider a single configuration ($\nTrain = 1$)
characterized by $\parami{i}=0$, $i=1,\ldots,9$ with no applied forcing
$\parami{i}=-2$, $i=9,\ldots,16$.  The nominal forces that affect the initial condition 
(see Table \ref{tab:parameterization}) are set to $\fNomnum{1} = \fNomnum{2} = 2\mathrm{kg} \times 9.81
\mathrm{m/s}^2$ and $\fNomnum{3} = \fNomnum{4} = 0.4 \mathrm{kg} \times 9.81
\mathrm{m/s}^2$. 
The time-step size is set to $\timestep=0.25$ seconds, leading to
$\totalTS=100$ total time steps.  This value was determined by a
timestep-verification study using a timestep-refinement factor of two; a
timestep of $0.25$ seconds led to an approximated rate of convergence in the
	output quantity $d$ at the end of the time interval of
	$1.40$ (which is reasonably close to the scheme's asymptotic rate of
	convergence of 2.0) and an approximated error in this quantity (computed via
	Richardson extrapolation) of $0.99\%$.

	We assess the performance of the ROMs with spatial-complexity reduction
	(i.e., Gappy POD and the structure-preserving ROM) using two different
	sets of sample indices. First, we set the number of sample nodes equal to
	20\% of the total nodes in the mesh (i.e., $n_s = 200$), which leads to
	$\nsample = 600$. We also employ a sampling fraction of 5\%, which leads to
	$\nsample = 150$.  	% ROM data:
	% gappyInfo_nominal_FREElow2_name_tr100_sz250_ts100_dt250tol3_p1_d15_maxFr1_nW10Galfore9new0_sn0IC_nR10_nIntF-20_nOff1_aj1.mat
	For the forecasting technique, the Newton-step
	threshold is set to $\newtonThreshold = 0$ and the maximum memory is set to
	$\memoryMax = 9$. For the `no forecasting' case, we employ a zeroth-order
	polynomial extrapolation. For experiments in this section, we declare the Newton iterations
	to have converged when the residual norm reaches $10^{-4}$ of its value
	computed using a zero acceleration and the values of the displacement and
	velocity at the beginning of the timestep.

\begin{figure}[htbp] \centering
	\subfigure[sampling fraction $\nsample/\ndof = 0.05$]{\label{fig:limiting5}
		\includegraphics[width=0.4\textwidth]{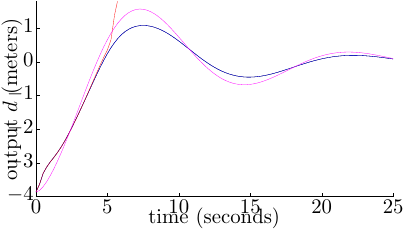}}
	\subfigure[sampling fraction $\nsample/\ndof = 0.20$]{
		\includegraphics[width=0.4\textwidth]{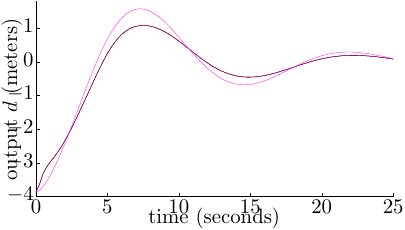}}
	\caption{Ideal case: Online responses for the full-order model (black,
	hidden), Galerkin ROM (blue) and Gappy POD ROM (red), and
	structure-preserving ROM (magenta) for different sampling fractions. Note
	that the Gappy POD ROM is unstable for a sampling fraction of 0.05.}\label{fig:limiting}
	\end{figure}

	The full-order-model simulation consumed $16.8$ minutes and incurred $229$
	Newton iterations ($\newtonItsAvg_\FOM = 2.29$). Table \ref{tab:limiting} and Figure \ref{fig:limiting}
	report the results for the reduced-order models.
	First, note that the relative errors generated by Galerkin and Gappy POD ROM
	with 20\% sampling are essentially zero.
	This is expected, because the reduced bases are not truncated and the inputs
	are fixed. Further, note that the Galerkin ROM without forecasting generates no speedup; this is
	expected because it is not equipped with a spatial-complexity-reduction technique (see
	Section \ref{sec:hyper}). The other two techniques---which employ
	spatial-complexity-reduction approximations---lead to speedups. The
	exception is Gappy POD with 5\% sampling, which yields an unstable
	response; this is depicted in Figure \ref{fig:limiting5}. For this reason,
	subsequent experiments employ a larger sampling fraction for the Gappy
	POD ROM compared with the structure-preserving ROM. 
		
	Importantly, note that the reduced-order models exhibit very little
	temporal-complexity reduction (i.e., $\temporalRedFactor \approx 1.0$) in
	the absence of the proposed forecasting technique.
	When the models employ the proposed forecasting technique, the number of
	Newton iterations decreases, leading to temporal-complexity
	reductions of $\temporalRedFactor=114.5$ for the Galerkin ROM and
	$\temporalRedFactor=2.26$  and $\temporalRedFactor=2.25$ for the
	best-performing Gappy POD and structure-preserving ROMs, respectively. In turn, this leads to 
	improved wall-time speedups in all cases.
	
The Galerkin ROM case presented here can be viewed as the best
possible performance for the method (applied to this problem): the temporal
behavior of the system is exactly predictable, as the inputs have not changed,
and the reduced basis has not been truncated. So, the forecast is `perfect'
after only one time step for the Galerkin ROM.
This means that for each time step after the first one, the initial guess
generated by the forecasting method is equal to the solution at that time
step, so no Newton steps are needed to compute the solution. As a result, no
Newton iterations are carried out beyond the first time step. The next
sections investigate the forecasting method's performance in the (more
realistic) case of varying inputs and truncated bases.

\begin{remark}
Note that the speedup (2.00) of the
	Galerkin ROM with forecating is not nearly as significant as the reduction
	factor (114.5), as Newton
	iterations are not the only aspect of the simulation that contribute to
	computational time. For example, the solution, velocity, and
	acceleration are updated at each time step, the residual is computed at each
	time step to check for convergence, outputs are computed, etc. We expect
	these two values to align more closely for problems where the computational
	cost of the Newton iterations dominates the overall simulation time.
\end{remark}
	 \begin{table}[htd] 
	 \centering 
	 \scriptsize{\begin{tabular}{|c|c||c|c|c|c||c|c|c|c|} 
		\hline 
	\multirow{4}{*}{ROM method}&\multirow{4}{1.5cm}{\centering{sampling
	fraction $\nsample/\ndof$}}&
	\multirow{4}{1.0cm}{\centering{relative error $\efomone$}}& 
	\multicolumn{3}{c||}{No forecasting} & \multicolumn{3}{|c|}{With forecasting} \\
	\cline{4-9}&&&
	\multirow{3}{1.0cm}{\centering{Newton its $\newtonItsAvg\totalTS$}}&
	\multirow{3}{1.0cm}{\centering{speedup $\speedup$}} &
	\multirow{3}{1.0cm}{\centering{reduction factor $\temporalRedFactor$}}&
	%\multirow{3}{1.0cm}{\centering{relative error $\efomone$}}& 
	 \multirow{3}{1.0cm}{\centering{Newton its $\newtonItsAvg\totalTS$}} & \multirow{3}{1.2cm}{\centering{speedup $\speedup$}} & \multirow{3}{1.0cm}{\centering{reduction factor $\temporalRedFactor$}}\\
	&&&& & & & &  \\
	&&&&& & & &  \\
		\hline 
	Galerkin &- & $8.93\times 10^{-6}$& 209 & 0.955&1.10& %$1.49\times 10^{-8}$&
	%$1.49\times 10^{-8}$&
	2 &2.00 &114.5 \\
	\hline
	\multirow{2}{*}{Gappy POD} &0.2& $1.60\times 10^{-5}$&209&2.97& 1.10&%$1.73\times 10^{-4}$  &
	%$1.73\times 10^{-4}$&
	101 &4.34 &2.26 \\
	& 0.05& unstable &-&-&- &%$1.73\times 10^{-4}$  &
	%$5.24\times 10^{-1}$&
	- &- &- \\
	\hline
	\multirow{2}{*}{structure preserving}& 0.2& $5.06\times 10^{-2}$&199 &3.40 &1.15 &% $2.81\times 10^{-1}$&
	%$2.81\times 10^{-1}$&
	107 &4.27 &2.14 \\
	& 0.05& $4.98\times 10^{-2}$&199 &12.7 &1.15 &% $2.81\times 10^{-1}$&
	%$2.77\times 10^{-1}$&
	102 &16.3 &2.25 \\
	\hline
		\end{tabular} }
		\caption{Ideal case: forecast performance.} 
		\label{tab:limiting} 
		%tab:~/projects/temporalComplexity/temporal_experiments_2_maxiter_redo/size_250/performanceText/gap_sparsenominal_FREElow2_nominal_FREElow2_name_tr100_sz250_ts100_dt250tol4_p1_d15_maxFr1_nW10Galfore9new0_sn0IC_nR10_nIntF-20_nOff1_aj1lbfFast1e_68nopiece_perform.txt
		% errors:
		% allnominal_FREElow2_nominal_FREElow2_name_tr100_sz250_ts100_dt250tol4_p1_d15_maxFr1_nW10Galfore0new0_sn0IC_nR10_nIntF-20_nOff1_aj1lbfFast1e_68nopiece_error.txt
			% 8.929151e-06     1.595056e-05    8.263829e-06    5.061822e-02
		% allnominal_FREElow2_nominal_FREElow2_name_tr100_sz250_ts100_dt250tol4_p1_d15_maxFr1_nW10Galfore0new0_sn0IC_nR10_nIntF-5_nOff1_aj1lbfFast1e_68nopiece_error.txt
			% 8.929151e-06     3.450889e-01    2.609366e-02    4.981478e-02
		\end{table} 
	 
	\subsection{Unforced, varying inputs}\label{sec:FREElow2}
\begin{figure}[htbp] 
		\centering 
	\includegraphics[width=0.4\textwidth]{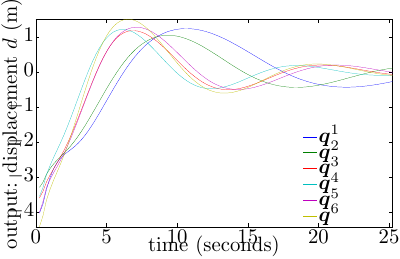}
			\caption{Unforced, varying inputs: Full-order model responses at training points in
			parameter space.} 
			 \label{fig:consPredict} 
				\end{figure}

	% Problem description
	%problem_numbers = [25];
	%dt = 0.25;
	%basisDim = 5;
	%basisDimRes = 9;
	%nIntF = -5;
	%nIntF = -60;
		%The  number of rows in $\restrict$ is set to be 60\% of the dofs
We now consider a fully predictive scenario with
$\paramOnline\not\in \paramTrain$. Again, we set the forces to zero, which
implies $\parami{i}=-2$, $i=9,\ldots,16$. We use $\nTrain = 6$ training points
and determine $\paramTrain$ using Latin hypercube sampling
\citep{mckay1979comparison}. We randomly select two online points.  Figure
\ref{fig:consPredict} depicts the tip displacement for the training points. As
the problem setup is the same as the previous section (except for the
parameter variation), we employ the same timestep size of
$\timestep=0.25$ seconds, leading to $\totalTS = 100$ time steps.

\begin{figure}[htbp] 
	\centering 
	\subfigure[online point $\paramOnlineNum{1}$]{
	\includegraphics[width=0.4\textwidth]{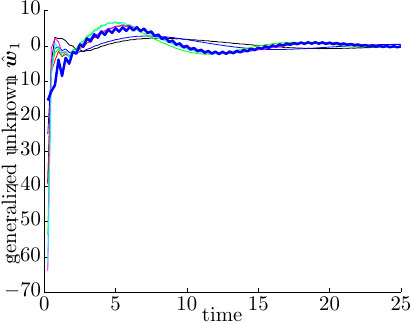} }
	\subfigure[online point $\paramOnlineNum{2}$]{
	\includegraphics[width=0.4\textwidth]{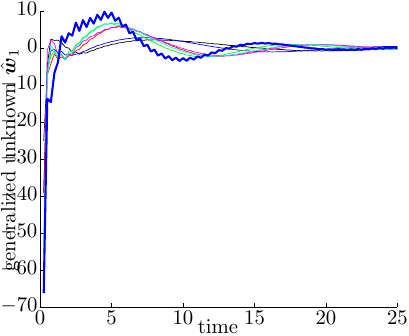} }
	\caption{Unforced, varying inputs: First generalized unknown at online point (bold
curve) and training points (thin curves).}
	\label{fig:firstGenCoordFREElow2} 
	%/Users/ktcarlb/projects/temporalComplexity/temporal_experiments_2_maxiter_redo/size_250/performanceText: 
	%%compare(3,1,0)
	%%compare(3,3,0)
	\end{figure}

To gain insight into the proposed method's potential, Figure
\ref{fig:firstGenCoordFREElow2} depicts the time evolution of the first generalized
unknown $\redUnknown_1$---which is one of the forecasted variables---for the online and training points. Importantly, note that the qualitative
response of this unknown is quite similar across parameter variation, which
suggests that the forecasting method has the potential to generate accurate
forecasts.

To construct the reduced-order models, we employ truncation critera of
$\nu_\state = 1-10^{-5}$, which leads to a basis dimension of $\nstate =8$, and
$\nu_{\residual} = 1-10^{-9}$ for Gappy POD, which results in a dimension
of 316 for the residual basis.  For the structure-preserving ROM, we sample 5\% of the indices such
that $\nsample = 150$; as this led to instabilities for Gappy POD, we sample
60\% of the indices (i.e., $\nsample = 1800$) for that method.
% gappyInfo_train6_FREElow2_tr100_sz250_ts100_dt250tol6_p1_d15_maxFr1_nW5Galfore0new0_sn0IC_nR9_nIntF-5_nOff1_aj1.mat

\begin{figure}[htbp] 
	\centering 
	\subfigure[online point $\paramOnlineNum{1}$]{
	\includegraphics[width=0.4\textwidth]{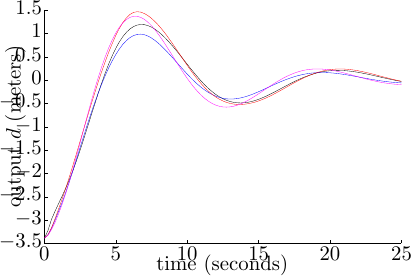} }
	\subfigure[online point $\paramOnlineNum{2}$]{
	\includegraphics[width=0.4\textwidth]{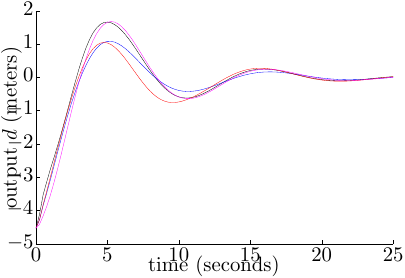} }
	\caption{Unforced, varying inputs: Online responses for the full-order model (black), Galerkin
	ROM (blue), Gappy POD ROM (red), and structure-preserving ROM
	(magenta).}
	\label{fig:FREElow2responses} 
	%/Users/ktcarlb/projects/temporalComplexity/temporal_experiments_2_maxiter_redo/size_250/matfigs.d
	\end{figure}

Figure \ref{fig:FREElow2responses} reports the responses of the full-order
model and all three reduced-order models. The full-order-model simulation
required 18.5 minutes and 307 total Newton iterations ($\newtonItsAvg_\FOM =
3.07$) for online point
%FOMrand1_FREElow2_sz250_ts100_dt250tol6_p1_d15_maxFr1.mat
$\paramOnlineNum{1}$ and 20.4 minutes and 347 Newton iterations ($\newtonItsAvg_\FOM =
3.47$) for online point
$\paramOnlineNum{2}$.
%FOMrand3_FREElow2_sz250_ts100_dt250tol6_p1_d15_maxFr1.mat
Note that the reduced-order models
are very accurate at the prediction points. At online point
$\paramOnlineNum{1}$, they generate relative errors
$\efomone$ of $3.33\times 10^{-2}$ (Galerkin), $2.56\times 10^{-2}$ (Gappy POD), and
$4.66\times 10^{-2}$ (structure preserving). At
online point $\paramOnlineNum{2}$, the relative errors are $3.48\times
10^{-2}$ (Galerkin), $4.07\times 10^{-2}$
(Gappy POD), and $2.45\times 10^{-2}$ (structure
preserving).\footnote{Different initial guesses for the Newton solver lead to
(slightly) different computed responses. Thus, the 
ROM responses in principle depend on the forecasting method. However, the
resulting differences in errors were negligible in these experiments;
therefore, we only report the ROM error generated by an initial guess of
zero (i.e., polynomial forecast, $\memory=0$ in Figure
\ref{fig:FREElow2responses}).}
% allrand1_FREElow2_train6_FREElow2_tr100_sz250_ts100_dt250tol6_p1_d15_maxFr1_nW5Galfore0new0_sn0IC_nR9_nIntF-5_nOff1_aj1lbfFast1e_68nopiece_error.txt
	% 3.329038e-02     4.574401e+01    1.380385e+00    4.657194e-02
% allrand1_FREElow2_train6_FREElow2_tr100_sz250_ts100_dt250tol6_p1_d15_maxFr1_nW5Galfore0new0_sn0IC_nR9_nIntF-60_nOff1_aj1lbfFast1e_68nopiece_error.txt
	% 3.329038e-02     2.556691e-02    2.429268e+02    4.666150e-02
% allrand3_FREElow2_train6_FREElow2_tr100_sz250_ts100_dt250tol6_p1_d15_maxFr1_nW5Galfore0new0_sn0IC_nR9_nIntF-5_nOff1_aj1lbfFast1e_68nopiece_error.txt
	% 3.476894e-02     3.341230e+01    2.279954e+01    2.453370e-02
% allrand3_FREElow2_train6_FREElow2_tr100_sz250_ts100_dt250tol6_p1_d15_maxFr1_nW5Galfore0new0_sn0IC_nR9_nIntF-60_nOff1_aj1lbfFast1e_68nopiece_error.txt
	% 3.476894e-02     4.067162e-02    1.646259e+02    2.461992e-02

%sparserand1_FREElow2_train6_FREElow2_tr100_sz250_ts100_dt250tol6_p1_d15_maxFr1_nW5Galfore0new0_sn0IC_nR9_nIntF-5_nOff1_aj1lbfFast1e_68.txt
%sparserand3_FREElow2_train6_FREElow2_tr100_sz250_ts100_dt250tol6_p1_d15_maxFr1_nW5Galfore0new0_sn0IC_nR9_nIntF-5_nOff1_aj1lbfFast1e_68.txt
%gappyrand1_FREElow2_train6_FREElow2_tr100_sz250_ts100_dt250tol6_p1_d15_maxFr1_nW5Galfore0new0_sn0IC_nR9_nIntF-60_nOff1_aj1.txt
%gappyrand3_FREElow2_train6_FREElow2_tr100_sz250_ts100_dt250tol6_p1_d15_maxFr1_nW5Galfore0new0_sn0IC_nR9_nIntF-60_nOff1_aj1.txt

\begin{figure}[htbp] 
	\centering 
	\subfigure[online point $\paramOnlineNum{1}$]{
	\includegraphics[width=0.8\textwidth]{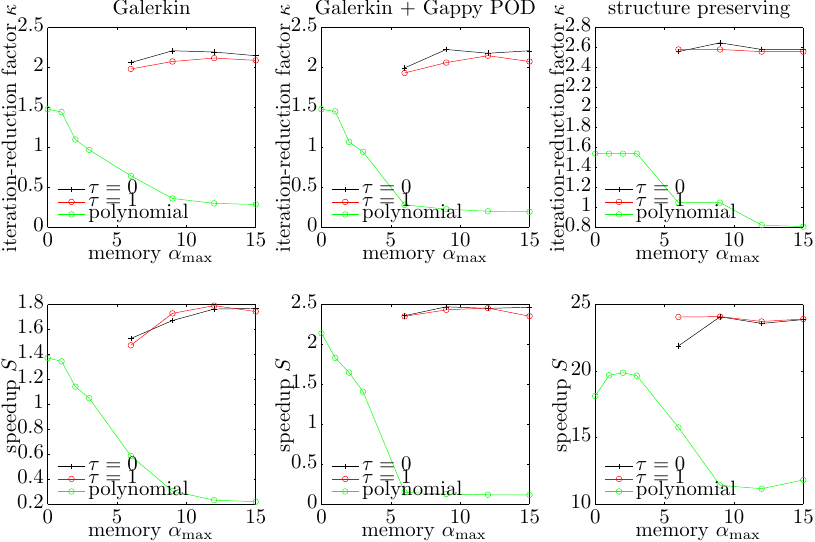} }
	\subfigure[online point $\paramOnlineNum{2}$]{
	\includegraphics[width=0.8\textwidth]{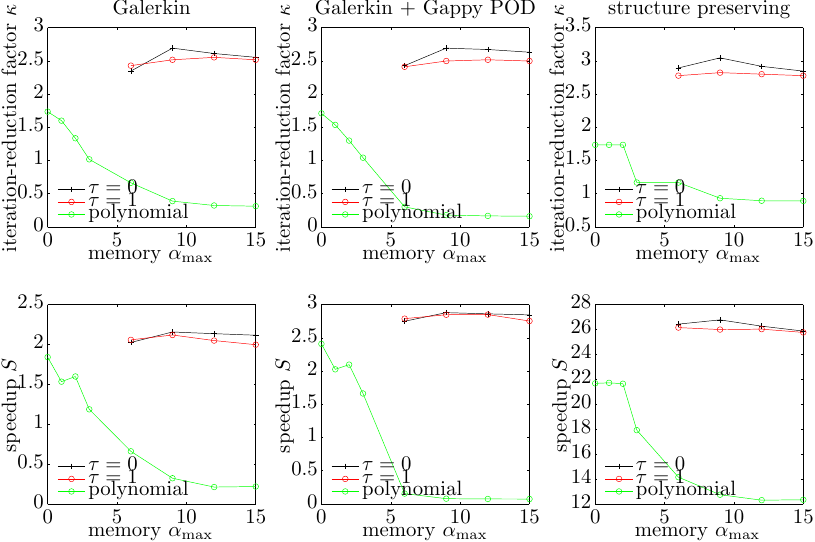} }
	\caption{Unforced, varying inputs: Performance of the forecasting method.
	The proposed forecasting method decreases both
	the number of requried Newton iterations and simulation time compared with
	polynomial extrapolation in nearly all cases.}
	\label{fig:FREElowFore} 
	%/Users/ktcarlb/projects/temporalComplexity/temporal_experiments_2_maxiter_redo/size_250/performanceText: 
	%%compare(5,1,0)
	%%compare(5,3,0)
	\end{figure}

Figure \ref{fig:FREElowFore} reports the Newton-iteration and wall-time
performance of the reduced-order models for different forecasting strategies
at the two online points.
First, note that the proposed forecasting method always yields better
performance than polynomial extrapolation, regardless of the values for the
forecasting parameters or polynomial degree. Second, observe that the
performance of the proposed forecasting method is relatively insensitive to
its parameters $\newtonThreshold$ and $\memoryMax$. Also, note that adding
`memory' to the polynomial extrapolation forecast---which yields a higher-degree
extrapolant---is almost always deleterious to its performance. In addition,
improvement in wall-time speedup provided by the forecasting technique is not as
strong as the improvement in number of Newton iterations. This can be
attributed to the presence of other operations (e.g., solution updating,
residual computation to check for convergence) that contribute to the
simulation time (see the remark in Section \ref{sec:limiting}). Finally,
observe that the speedups generated by the structure-preserving method are far
superior to those generated by Galerkin and Gappy POD. This is due to the fact
that the structure-preserving method employed only $\nsample = 150$, whereas
Galerkin is not equipped with a spatial-complexity-reduction mechanism and
Gappy POD required $\nsample = 1800$ to generate a stable response.

	\subsection{Forced, varying inputs}\label{sec:FRCSTRlow2}
\begin{figure}[htbp] 
  \centering 
	\subfigure[Section \ref{sec:FRCSTRlow2}]{
		 \label{fig:forcedFOMtraininglow} 
\includegraphics[width=0.4\textwidth]{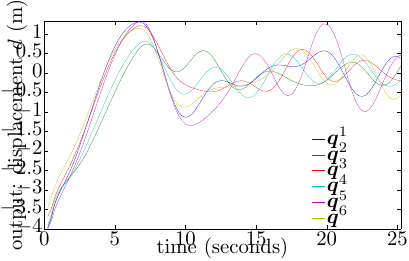}
}
	\subfigure[Section \ref{sec:FRCSTR2}: larger parameter variation]{
		 \label{fig:forcedFOMtraining} 
\includegraphics[width=0.4\textwidth]{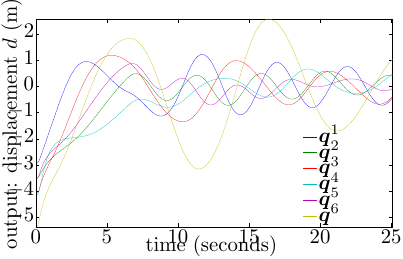}
}
	  \caption{Forced, varying inputs: Full-order model responses at training points in
		parameter space. Note that larger parameter variation leads to larger
		parameter-induced changes in the output.} 
		  \end{figure}

%FRCSTRlow2
	In this section, we activate the external forcing, thereby allowing
	$\parami{i}\in\left[-0.5,0.5\right]$, $i=1,\ldots,16$. The timestep was again
	set to $\timestep = 0.25$ seconds, leading to $\totalTS = 100$ time steps.
	This value was again determined by a timestep-verification study at the
	nominal configuration $\paramNom$ using a refinement factor of two. The
	approximated rate of convergence in the output quantity at the end of the
	time interval for this timestep size was determined to be 1.67 (close to the
	asymptotic value of 2.0), and the
	error in this quantity as approximated by Richardson extrapolation was
	$1.33\%$. As before, we used Latin hypercube sampling to determine the
	$\nTrain=6$ training points; Figure \ref{fig:forcedFOMtraininglow} reports the full-order-model responses
at these points. Note that parameter variation leads to significant changes in
the response. We randomly select two online points at which we will perform
prediction with the ROMs.
% Parameter variation

\begin{figure}[htbp] 
\centering 
\subfigure[online point $\paramOnlineNum{1}$]{
\includegraphics[width=0.4\textwidth]{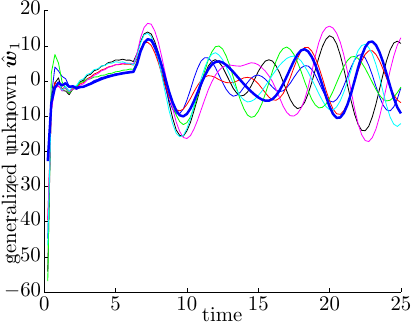} }
\subfigure[online point $\paramOnlineNum{2}$]{
\includegraphics[width=0.4\textwidth]{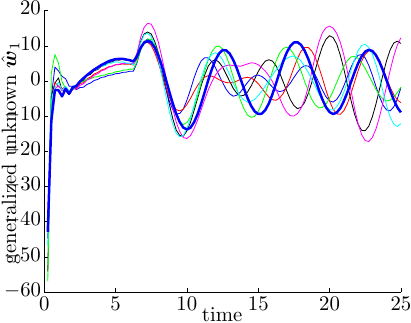} }
\caption{Forced, varying inputs: First generalized unknown at online point (bold
curve) and training points (thin curves).}
\label{fig:firstGenCoordFRCSTRlow2} 
%/Users/ktcarlb/projects/temporalComplexity/temporal_experiments_2_maxiter_redo/size_250/figures.d/coordPlot/: 
%%compare(3,1,0)
%%compare(3,3,0)
\end{figure}

Figure \ref{fig:firstGenCoordFRCSTRlow2} depicts the time evolution of the first generalized
unknown $\redUnknown_1$ for the online and training points. As before, there
is qualitative similarity of this forecasted variable for the different
points; this suggests the forecasting method can again realize computational
savings. Also, note that the character of the response changes appreciably
when the external force is activated at $t = 6.25$ seconds.

The reduced-order models employ trunction critera of $\nu_\state = 1-10^{-6}$
(basis dimension of $\nstate = 16$) and a residual-basis dimension of 1800. The
structure-preserving method approximates the external force via Gappy POD (see
Ref.~\citep{carlberg2012spd}); for this purpose, it employs a truncation
criterion of $\nu_{\fext} = 1$, leading to a basis dimension of
4.\footnote{Note that the external force is composed of only four linearly
independent components $\forceDist_i$, $i=1,\ldots, 4$ (see
Eq.~\eqref{eq:externalForceOne}).} Again, the Gappy POD ROM employs a sampling
rate of $60\%$ ($\nsample = 1800$) and the structure-preserving
ROM employs a samping percentage of $5\%$ ($\nsample = 150$).\footnote{The Gappy POD
ROM was unstable for $\nsample = 150$.}

\begin{figure}[htbp] 
\centering 
\subfigure[online point $\paramOnlineNum{1}$]{
\includegraphics[width=0.4\textwidth]{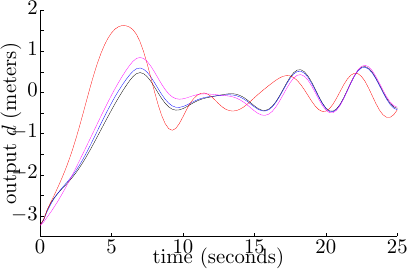} }
\subfigure[online point $\paramOnlineNum{2}$]{
\includegraphics[width=0.4\textwidth]{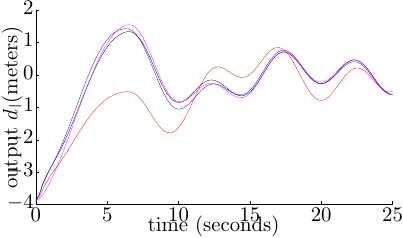} }
\caption{Forced, varying inputs: Online responses for the full-order model (black), Galerkin
ROM (blue), Gappy POD ROM (red), and structure-preserving ROM
(magenta).}
\label{fig:FRCSTRlow2onlineROM} 
%/Users/ktcarlb/projects/temporalComplexity/temporal_experiments_2_maxiter_redo/size_250/matfigs.d
\end{figure}

Figure \ref{fig:FRCSTRlow2onlineROM} reports the responses of the full-order
model and the reduced-order models at the online prediction points. The
full-order model consumed 20.3 minutes and 330 Newton iterations ($\newtonItsAvg_\FOM =
3.3$) at online point
%FOMrand2_FRCSTRlow2_sz250_ts100_dt250tol6_p1_d15_maxFr3.mat
$\paramOnlineNum{1}$ and 22.6 minutes and 360 Newton iterations ($\newtonItsAvg_\FOM =
3.6$) at point
%FOMrand4_FRCSTRlow2_sz250_ts100_dt250tol6_p1_d15_maxFr3.mat
$\paramOnlineNum{2}$. The
relative errors $\efomone$ of the ROMs at online point $\paramOnlineNum{1}$
are $1.56\times 10^{-2}$ (Galerkin), $1.56\times 10^{-1}$ (Gappy POD), and
$5.78\times 10^{-2}$ (structure-preserving). For
online point $\paramOnlineNum{2}$, the errors are $2.41\times 10^{-2}$
(Galerkin), $1.68\times 10^{-1}$ (Gappy
POD), and $3.05\times 10^{-2}$ (structure-preserving). Note that the Galerkin and
structure-preserving ROMs are quite accurate, but the Gappy POD ROM incurs
significant errors. 
% allrand2_FRCSTRlow2_train6_FRCSTRlow2_tr100_sz250_ts100_dt250tol6_p1_d15_maxFr3_nW6Galfore0new0_sn0IC_nR10_nIntF-5_nOff1_aj1lbfFast1e_68nopiece_error.txt
	% 1.561335e-02     1.822106e+00    2.029367e+00    5.784881e-02
% allrand2_FRCSTRlow2_train6_FRCSTRlow2_tr100_sz250_ts100_dt250tol6_p1_d15_maxFr3_nW6Galfore0new0_sn0IC_nR10_nIntF-60_nOff1_aj1lbfFast1e_68nopiece_error.txt
	% 1.561335e-02     1.558101e-01    9.517481e-01    4.845890e-02
% allrand4_FRCSTRlow2_train6_FRCSTRlow2_tr100_sz250_ts100_dt250tol6_p1_d15_maxFr3_nW6Galfore0new0_sn0IC_nR10_nIntF-5_nOff1_aj1lbfFast1e_68nopiece_error.txt
	%2.414938e-02     1.700117e-01    1.801071e+00    3.047891e-02
% allrand4_FRCSTRlow2_train6_FRCSTRlow2_tr100_sz250_ts100_dt250tol6_p1_d15_maxFr3_nW6Galfore0new0_sn0IC_nR10_nIntF-60_nOff1_aj1lbfFast1e_68nopiece_error.txt
	% 2.414938e-02     1.676560e-01    8.229441e-01    2.109468e-02

\begin{figure}[htbp] 
\centering 
\subfigure[online point $\paramOnlineNum{1}$]{
\includegraphics[width=0.8\textwidth]{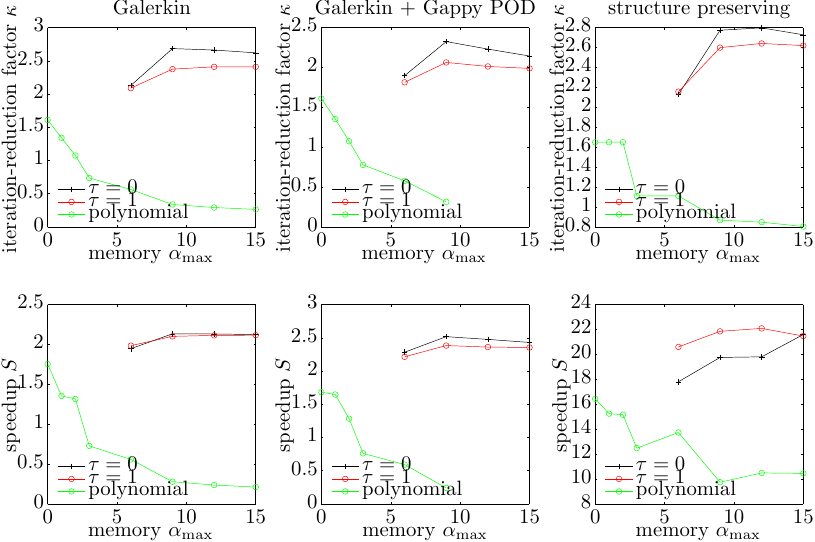} }
\subfigure[online point $\paramOnlineNum{2}$]{
\includegraphics[width=0.8\textwidth]{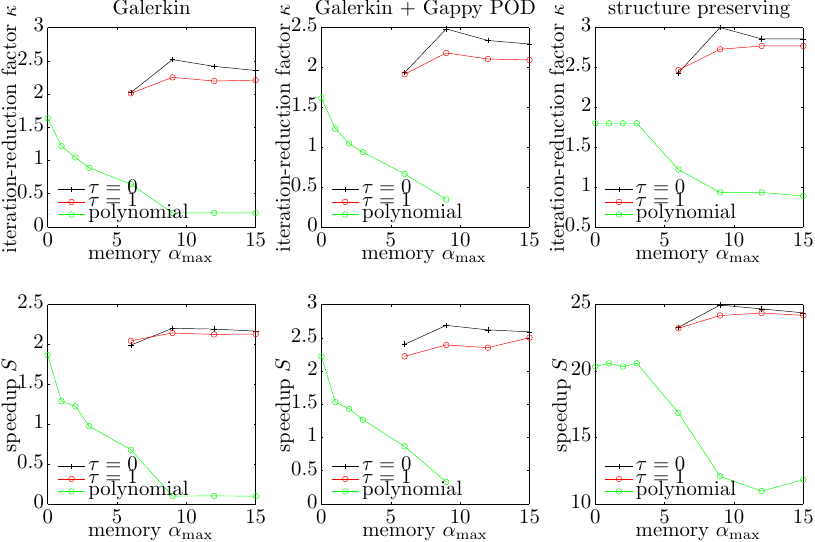} }
\caption{Forced, varying inputs: Performance of the forecasting method.
For all reduced-order models, the proposed forecasting method decreases both
the number of requried Newton iterations and simulation time compared with
polynomial extrapolation.}
\label{fig:FRCSTRlowFore} 
%/Users/ktcarlb/projects/temporalComplexity/temporal_experiments_2_maxiter_redo/size_250/performanceText:
%%compare(4,2,0)
%%compare(4,4,0)
\end{figure}

Figure \ref{fig:FRCSTRlowFore} reports the Newton-iteration and wall-time
performance of the ROMs for different forecasting strategies. The results are
very similar to those for the unforced case: the proposed forecasting method
nearly always exhibits performance superior to that of 
polynomial extrapolation, the proposed method is relatively insensitive to the parameters
$\newtonThreshold$ and $\memoryMax$, and high-order polynomial extrapolation performs
very poorly. In addition, improvement in
iteration-reduction factor $\temporalRedFactor$ exceeds the improvement in
speedup $\speedup$, and the
structure-preserving method generates the largest speedups due to the fact it
employs the smallest number of sample indices. Additionally, notice the
`missing' data points for polynomial extrapolation with $\memoryMax = 12$ and
$\memoryMax = 15$ for the Gappy POD ROM; these missing data indicate that the
Gappy POD ROM did not converge for these forecasts. This
implies that the initial guesses were so poor that the globalized Newton method
failed to generate an acceptable solution within the alloted 500 Newton
iterations at least one time step.

Also, note that employing $\newtonThreshold=0$ appears to systematically
outperform $\newtonThreshold = 1$ in terms of the iteration-reduction factor
$\temporalRedFactor$ metric. However, this does not always lead to an
improvement in speedup (see the structure-preserving ROM for
$\paramOnlineNum{1}$). This can be attributed to the fact that employing
$\newtonThreshold=0$ results in more frequent forecast recomputation (i.e.,
whenever the number of Newton iterations exceeds zero) than the
$\newtonThreshold = 1$ case.

These results highlight that the proposed forecasting method is applicable
even for the more challenging problem of parameterized forced responses.

\subsection{Forced, varying inputs, larger parameter variation.}\label{sec:FRCSTR2}
In this section, we assess the performance of the method for the same problem
as Section \ref{sec:FRCSTRlow2}, but with an increased parameter variation,
i.e., $\paramDomain = \left[-1,1\right]^{16}$. This poses a greater challenge
for both the reduced-order models and the forecasting method, as they now rely
on training data from the same number of points (we keep $\nTrain = 6$) in a
larger parameter domain. As the model now undergoes larger parameter
variation, we decrease the timestep size to $\timestep = 0.1$ seconds,
leading to $\totalTS = 250$ total time steps.\footnote{The full-order model
did not converge for several of the training points when $\timestep = 0.25$
seconds was employed.} Note that this timestep remains in the asymptotic
range of convergence for the nominal configuration $\paramNom$, as it is
smaller than the previously verified value of 0.25 seconds. Again, training
points are chosen by Latin hypercube sampling, and the online points are
selected randomly. Figure \ref{fig:forcedFOMtraining} reports the
full-order-model responses at the training points; note that the changes in
the response are in fact more significant than for the previous case with
smaller parameter variation.

\begin{figure}[htbp] 
\centering 
\subfigure[online point $\paramOnlineNum{1}$]{
\includegraphics[width=0.4\textwidth]{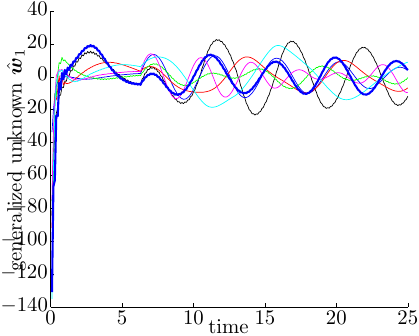} }
\subfigure[online point $\paramOnlineNum{2}$]{
\includegraphics[width=0.4\textwidth]{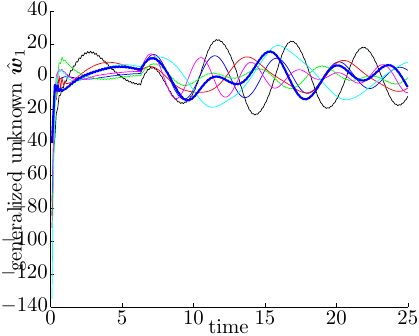} }
\caption{Forced, varying inputs, larger parameter variation: Parameter dependence of the first generalized coordinate.}
\label{fig:firstGenCoordFRCSTR2} 
%/Users/ktcarlb/projects/temporalComplexity/temporal_experiments_2_maxiter_redo/size_250/performanceText: 
%%compare(3,1,0)
%%compare(3,3,0)
\end{figure}

Figure \ref{fig:firstGenCoordFRCSTR2} again reports the time evolution of
the first generalized unknown. Note that again there is similar qualitative
structure across parameter variation. 
%By comparing with Figure
%\ref{fig:firstGenCoordFRCSTRlow2}, one can observe that the parameter-induced
%variation in the first generalized coordinate's response is larger in this
%case than when the parameter variation was smaller.

The attributes for the reduced-order models are the same as in Section
\ref{sec:FRCSTRlow2}, with one exception: a criterion of $\nu_\state =
1-10^{-4}$ is employed for the state, which associates with a  basis dimension
of $\nstate = 10$. Note that basis dimension is larger than in the previous
case.

\begin{figure}[htbp] 
\centering 
\subfigure[online point $\paramOnlineNum{1}$]{
\includegraphics[width=0.4\textwidth]{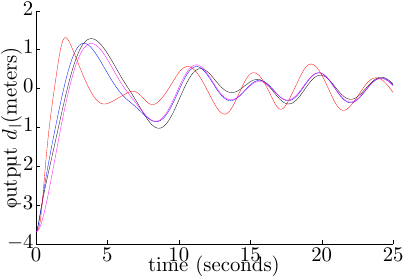} }
\subfigure[online point $\paramOnlineNum{2}$]{
\includegraphics[width=0.4\textwidth]{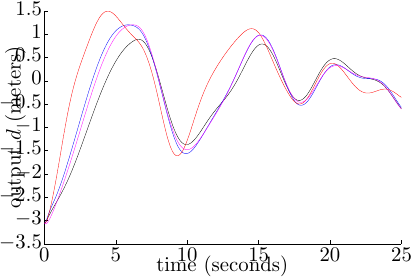} }
\caption{Forced, varying inputs, larger parameter variation: Online responses for the full-order model (black), Galerkin
ROM (blue), Gappy POD ROM (red), and structure-preserving ROM
(magenta)..}
\label{fig:FRCSTR2responses} 
%/Users/ktcarlb/projects/temporalComplexity/temporal_experiments_2_maxiter_redo/size_250/matfigs.d
\end{figure}

Figure \ref{fig:FRCSTR2responses} depicts the full-order-model response along
with those for the reduced-order models. 
The full-order model took 36.7 minutes and 605 Newton iterations ($\newtonItsAvg_\FOM =
2.42$) at online point
%FOMrand2_FRCSTR2_sz250_ts250_dt100tol6_p1_d15_maxFr3.mat
$\paramOnlineNum{1}$ and 37.3 minutes and 595 Newton iterations ($\newtonItsAvg_\FOM =
2.38$) at point
%FOMrand4_FRCSTR2_sz250_ts250_dt100tol6_p1_d15_maxFr3.mat
$\paramOnlineNum{2}$.
As before, the Galerkin and
structure-preserving ROMs are more accurate than the Gappy POD ROMs. The
relative errors $\efomone$ at point $\paramOnlineNum{1}$ are
$4.19\times 10^{-2}$ (Galerkin), $1.29\times 10^{-1}$ (Gappy POD), and
$3.65\times 10^{-2}$ (structure preserving). At online
point $\paramOnlineNum{2}$, the associated errors are 
$6.91\times 10^{-2}$ (Galerkin), $1.73\times 10^{-1}$ (Gappy POD), and
$5.67\times 10^{-2}$ (structure preserving).

%allrand2_FRCSTR2_train6_FRCSTR2_tr250_sz250_ts250_dt100tol6_p1_d15_maxFr3_nW4Galfore0new0_sn0IC_nR10_nIntF-5_nOff1_aj1lbfFast1e_68nopiece_error.txt
	%4.192957e-02     1.329337e+02    1.880790e+00    3.654973e-02
%allrand2_FRCSTR2_train6_FRCSTR2_tr250_sz250_ts250_dt100tol6_p1_d15_maxFr3_nW4Galfore0new0_sn0IC_nR10_nIntF-60_nOff1_aj1lbfFast1e_68nopiece_error.txt
	% 4.192957e-02     1.291437e-01    3.849880e+00    3.659511e-02
%allrand4_FRCSTR2_train6_FRCSTR2_tr250_sz250_ts250_dt100tol6_p1_d15_maxFr3_nW4Galfore0new0_sn0IC_nR10_nIntF-5_nOff1_aj1lbfFast1e_68nopiece_error.txt
	%6.910123e-02     3.312786e+02    5.638253e+00    5.667195e-02
%allrand4_FRCSTR2_train6_FRCSTR2_tr250_sz250_ts250_dt100tol6_p1_d15_maxFr3_nW4Galfore0new0_sn0IC_nR10_nIntF-60_nOff1_aj1lbfFast1e_68nopiece_error.txt
	%6.910123e-02     1.729902e-01    5.705433e+01    5.417196e-02

\begin{figure}[htbp] 
\centering 
\subfigure[online point $\paramOnlineNum{1}$]{
\includegraphics[width=0.8\textwidth]{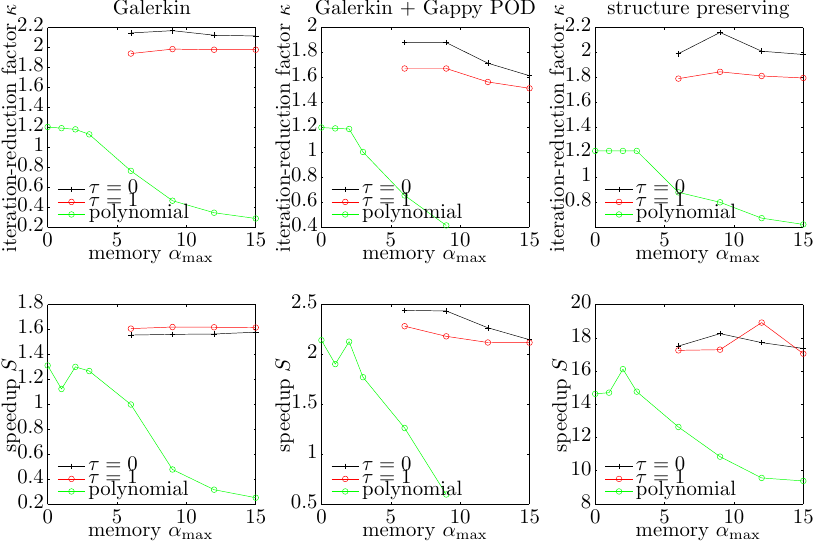} }
\subfigure[online point $\paramOnlineNum{2}$]{
\includegraphics[width=0.8\textwidth]{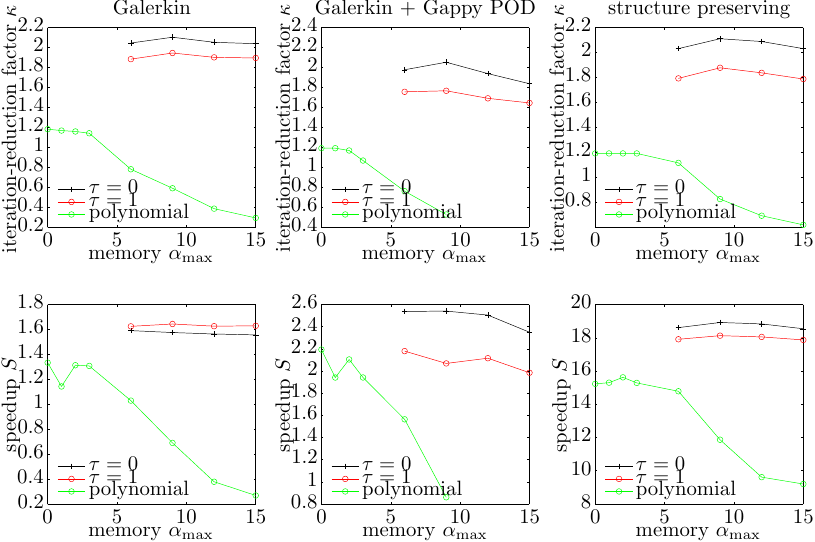} }
\caption{Forced, varying inputs, larger parameter variation: Performance of the forecasting method.
For all reduced-order models, the proposed forecasting method decreases both
the number of required Newton iterations and simulation time compared with
polynomial extrapolation.}
\label{fig:FRCSTRFore} 
%/Users/ktcarlb/projects/temporalComplexity/temporal_experiments_2_maxiter_redo/size_250/performanceText:
%%compare(2,2,0)
%%compare(2,4,0)
\end{figure}

Figure \ref{fig:FRCSTRFore} reports the Newton-iteration and wall-time results
for the different forecasting strategies. Note that the results are extremely
similar to those in Section \ref{sec:FRCSTRlow2}. The primary exception can
be seen by comparing Figure \ref{fig:FRCSTRFore} with
\ref{fig:FRCSTRlowFore}: the iteration-reduction factor $\temporalRedFactor$
and speedup $\speedup$ performance of the reduced-order models has
decreased. This can be attributed to the challenge of larger parameter
variation, as the ROMs are now responsible for capturing a wider range of
physics. 

From this set of experiments, we conclude that the proposed technique can
improve ROM performance even for problems with relatively large parameter
variation.

\subsection{Average performance}\label{sec:avg}

\begin{figure}[htbp] 
\centering 
\subfigure[Reduction-factor improvement]{
\includegraphics[width=0.4\textwidth]{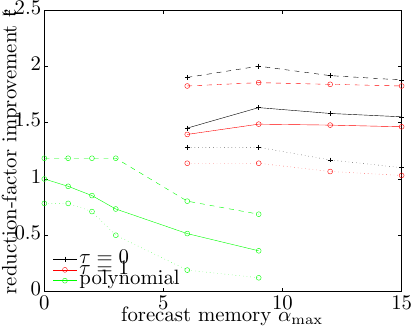} }
\subfigure[Speedup improvement]{
\includegraphics[width=0.4\textwidth]{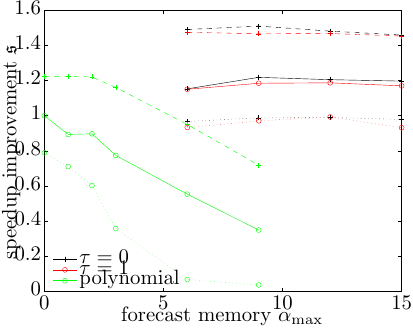} }
\caption{Performance of forecasting methods quantified over all reduced-order
models, problems, and online prediction points. The mean (solid line), maximum (dashed
line), and minimum (dotted) are reported.}
\label{fig:agglomerateResults} 
%/Users/ktcarlb/projects/temporalComplexity/temporal_experiments_2_maxiter_redo/size_250/performanceText:
%%agglomerateResults
\end{figure}

Finally, we summarize the performance of the forecasting techniques over the
complete set of experiments.
Figure \ref{fig:agglomerateResults} reports average, minimum, and maximum values of the
reduction-factor improvement $\temporalRedFactorImprove$, and speedup
improvement $\speedupImprove$ over all experiments (i.e., all three
experiments in Sections \ref{sec:FREElow2}--\ref{sec:FRCSTR2}, all three
reduced-order models, and both online points $\paramOnlineNum{1}$ and
$\paramOnlineNum{2}$).
Here, $\temporalRedFactorImprove =
\temporalRedFactor/\temporalRedFactor_{\mathrm{no}}$ and 
$\speedupImprove = \speedup/ \speedup_{\mathrm{no}}$ can each be computed for
a given ROM simulation; a subscript `no' indicates the value of the variable
for a zero initial guess (i.e., polynomial extrapolation with $\memory = 0$).
First, note that the proposed method always outperforms polynomial forecasting
in the mean, maximum, and minimum achieved performance for both
reduction-factor improvement $\temporalRedFactorImprove$ and speedup
improvement $\speedupImprove$.
Secondly, the maximum, minimum, and average
performance of polynomial forecasting were all made worse by increasing the
polynomial degree.

Finally, the best average performance was achieved for a forecast memory of
$\memoryMax = 9$ and Newton-iteration criterion of $\newtonThreshold = 0$. In
this case, the iteration-reduction factor was increased by $63\%$ on average;
the speedup was improved by $22\%$ on average.  Critically, note that these
temporal-complexity gains incur \emph{no additional error}, and so they
strictly serve to improve the performance of the ROMs with no penalty.

\section{Conclusions}

This paper has described a method for decreasing the temporal complexity of
nonlinear reduced-order models in the case of implicit time integration. The
method exploits knowledge of the dynamical system's temporal behavior in the
form of `time-evolution bases'; one such basis is generated for each generalized coordinate of the time
integrator's unknown during the (offline) training stage. During the (online)
deployed stage, these time-evolution bases are used---along with the solution at
recent time steps---to forecast the unknown at future time steps via Gappy
POD. If this forecast is accurate, the Newton-like solver will converge in
very few iterations, leading to computational-cost savings.

Numerical experiments demonstrated the potential of the method to
significantly improve the performance of nonlinear reduced-order models, even
in the presence of high-frequency content in the dynamics. The experiments
also demonstrated the effect of input parameters 
on the method's performance, and provided a parameter study to analyze the
effect of the method's parameters.

Future work includes devising a way to directly handle frequency and phase
shifts in the response, as well as time-shifted temporal behavior.
\section*{Acknowledgments}

%\marginpar{ktc - what other acknowledgments do we
%need?}
The authors acknowledge Julien Cortial for providing the original
nonlinear-truss code that was modified to generate the numerical results, as
well as the anonymous reviewers for their insightful suggestions.

This research was supported in part by an appointment to the Sandia National
Laboratories Truman Fellowship in National Security Science and Engineering,
sponsored by Sandia Corporation (a wholly owned subsidiary of Lockheed Martin
Corporation) as Operator of Sandia National Laboratories under its U.S.
Department of Energy Contract No. DE-AC04-94AL85000.

\appendix

\section{Implicit time-integration schemes: first-order ODEs}\label{app:implicitFirst}
For notational simplicity, consider a system without parametric inputs $\param$,
and define
 $\bar \f (\state,t)\equiv \f \left(\state;t,\forcing\left(t\right)\right)$ such that
 \begin{equation} 
 \dot \state = \bar \f \left(\state,t\right).
  \end{equation} 
 Further, denote by $h$ the time-step size at time step $n$.

\subsection{Implicit linear multi-step schemes}

A linear $k$-step method applied to first-order ODEs can be expressed as
 \begin{equation} 
 \sum_{j=0}^k\alpha_j\state^{n-j} = \timestep\sum_{j=0}^k\beta_j\bar \f \left(\state^{n-j},
 t^{n-j}\right),
 \end{equation} 
 where $\alpha_0\neq 0$ and $\sum\limits_{j=0}^k\alpha_j = 0$ is necessary for
 consistency. These methods are implicit if $\beta_0\neq 0$. In this case, the
 form of the 
 residual is
 \begin{equation} 
 \Rn\left(\unknown^n\right) = \alpha_0\unknown^n - \timestep\beta_0\bar \f (\unknown^n,t^n)+
 \sum_{j=1}^k\alpha_j\state^{n-j} -
 \timestep\sum_{j=1}^k\beta_j\bar \f \left(\state^{n-j},t^{n-j}\right)
  \end{equation} 
	and the explicit state update is simply
 \begin{equation} 
 \state^n = \unknown^n.
  \end{equation} 
Therefore, the unknown is the state at time $t^n$.

\subsection{Implicit Runge--Kutta schemes}\label{sec:rungeKutta}
For an $s$-stage Runge--Kutta scheme, the form of the residual is 
\begin{equation} 
\Rn_i\left(\unknownnone,\ldots,\unknownns\right) = \unknownni -
\bar \f (\state^{n-1} + \timestep\sum_{j=1}^sa_{ij}\unknownni, t^{n-1} +
c_ih),\quad i=1,\ldots, s
\end{equation} 
with the following explicit computation of the state:
\begin{equation} 
\state^{n} = \state^{n-1} + \timestep\sum_{i=1}^sb_i\unknownni.
\end{equation} 
 The unknowns correspond to the velocity $\dot \state$ at
 times $t^{n-1} + c_ih$, $i=1,\ldots, s$.

\section{Implicit time-integration schemes: second-order ODEs}
\label{app:implicitSecond}
For notational simplicity, consider a second-order differential equations
without parametric inputs $\param$ and define
$\bar \g\left(\state,\dot\state,t\right) \equiv
\g\left(\state,\dot\state;t,p(t)\right)$  such that
 \begin{equation} 
 \ddot \state = \bar \g\left(\state,\dot \state, t\right).
  \end{equation} 

% KTC: I could not find a good reference for implicit linear multi-step
% schemes applied to second-order differential equations where g depends on
% the velocity
%\subsection{Implicit linear multi-step schemes}
%
%A linear $k$-step method applied to second-order ODEs can be expressed as
% \begin{equation} 
% \sum_{j=0}^k\alpha_j\state^{n-j} = h^2\sum_{j=0}^k\beta_j\bar
% g\left(\state^{n-j},\dot\state^{n-j},
% t^{n-j}\right).
% \end{equation} 
% As before, $\alpha_0\neq 0$ and $\sum\limits_{j=0}^k\alpha_j = 0$ is necessary for
% consistency. These methods are implicit if $\beta_0\neq 0$. In this case, the form of the
% residual is (KTC FIX!)
% \begin{equation} 
% \Rn\left(\unknown^n\right) = \alpha_0\unknown^n - h^2\beta_0\bar g(\unknown^n,t^n)+
% \sum_{j=1}^k\alpha_j\state^{n-j} -
% h^2\sum_{j=1}^k\beta_j\bar g\left(\state^{n-j},t^{n-j}\right)
%  \end{equation} 
%	and the explicit state update is simply
% \begin{equation} 
% \state^n = \unknown^n.
%  \end{equation} 
%The unknown is the state at time $t^n$.
  
\subsection{Implicit Nystr\"om method}
Nystr\"om methods are partitioned Runge--Kutta schemes applied to second-order
ODEs. They lead to the following representation for the residuals:
\begin{align} \label{eq:nystrom}
\begin{split}
\Rn_i\left(\unknownnone,\ldots,\unknownns\right) =& \unknownni -\\
&\bar \g\left(\state^{n-1} + c_i\timestep\dot \state^{n-1}+h^2 \sum_{j=1}^s\bar
a_{ij}\unknownni, \dot \state^{n-1} + \timestep\sum_{j=1}^s\hat
a_{ij}\unknownni,t^{n-1}+ c_i\timestep\right),
\end{split}
\end{align} 
$i=1,\ldots,s$. The state and velocity are updated explicitly as
 \begin{gather} 
 \label{eq:nystromState}\state^n = \state^{n-1} + \timestep\dot\state^{n-1} + h^2\sum_{i=1}^s\bar b_i\unknownni\\
\label{eq:nystromVel} \dot\state^n = \dot\state^{n-1} + \timestep\sum_{i=1}^s\hat b_i\unknownni.
  \end{gather} 
	The unknowns correspond to the acceleration $\ddot \state$ at 
 times $t^{n-1} + c_ih$, $i=1,\ldots, s$.

\subsection{Implicit Newmark method}

The implicit Newmark method leads to the following residuals:
\begin{align}
\Rn(\unknown^n) = \unknown^n - \bar \g\left(\state^{n-1}+\timestep\dot \state^{n-1} +
\frac{h^2}{2}\left[\left(1 - 2\beta\right)\ddot \state^{n-1} +
2\beta\unknown^n\right],\dot \state^{n-1} +
	\timestep\left[\left(1-\gamma\right)\ddot\state^{n-1} + \gamma \unknown^n\right],t^n\right)
	\end{align}
	The state and velocity are explicitly updated as
	 \begin{gather} 
 \state^n = \state^{n-1}+\timestep\dot \state^{n-1} +
\frac{h^2}{2}\left[\left(1 - 2\beta\right)\ddot \state^{n-1} +
2\beta\unknown^n\right]\\
\dot\state^n = \dot \state^{n-1} +
\timestep\left[\left(1-\gamma\right)\ddot\state^{n-1} + \gamma \unknown^n\right].
\end{gather} 
Here, the unknown corresponds to the acceleration $\ddot \state$ at time $t^n$.

\section{Proper orthogonal decomposition} \label{app:POD}
Algorithm \ref{PODSVD} describes the method for computing a
proper-orthogonal-decomposition (POD) basis given a set of snapshots. The method
essentially amounts to computing the singular value decomposition of the
snapshot matrix. The left singular vectors define the POD basis.

\begin{algorithm}[htbp]
\caption{Proper-orthogonal-decomposition basis computation (normalized
snapshots)}
\begin{algorithmic}[1]\label{PODSVD}
\REQUIRE Set of snapshots $\snapsNo\equiv\{\w _i\}_{i=1}^\nsnap\subset\RR{\ndof }$,
energy criterion $\energyCrit\in[0,1]$
\ENSURE  $\podArgs{\snapsNo}{\energyCrit}$
\STATE\label{step:SVD} Compute thin singular value decomposition $
\snapmat= \boldsymbol U \boldsymbol \Sigma \boldsymbol V^T $, where
$\boldsymbol W\equiv\left[\w _1/\|\w _1\|\ \cdots\
\w _{n_\w }/\|\w _{n_\w }\|\right]$.
\STATE Choose dimension of truncated basis 
$\nstate = \nenergy(\nu)$, where 
 \begin{align} 
 \nenergy(\nu) &\equiv \arg\min_{i\in \mathcal V(\nu)}i\\
 \mathcal V(\nu)&\equiv \{n\in\{1,\ldots,\nsnap\}\ | \
 \sum_{i=1}^n\sigma_i^2/\sum_{i=1}^{\nsnap}\geq\nu\},
  \end{align} 
	and $\boldsymbol \Sigma \equiv \text{diag}\left(\sigma_i\right)$.
\STATE $\podArgs{\snapsNo}{\energyCrit}=\vecmat{\boldsymbol u}{\nstate}$,
where $\boldsymbol U
\equiv\vecmat{\boldsymbol u}{\nsnap}$.
\end{algorithmic}
\end{algorithm}
\bibliography{references,refs_jr}
\bibliographystyle{elsarticle/model1-num-names}
\end{document}